\documentclass[10 pt, journal]{IEEEtran}  
\IEEEoverridecommandlockouts

\usepackage{xspace,amssymb,epsfig,syntonly}
\usepackage{epsfig,amsmath,color}
\usepackage{xspace,syntonly,empheq}
\usepackage{url}
\usepackage{bbm}
\usepackage{mathtools}
\usepackage{cite}
\usepackage{graphicx}
\usepackage{caption}
\usepackage{subcaption}
\usepackage{algorithm,algorithmic}
\usepackage{verbatim}
\usepackage{amssymb}
\usepackage{amsthm}
\usepackage{circuitikz}
\usepackage{enumerate}

\newcommand{\utwi}[1]{\mbox{\boldmath $#1$}}

\newcommand{\diag}{\mathsf{diag}}

\renewcommand{\hat}{\widehat}
\renewcommand{\tilde}{\widetilde}

\newcommand{\cD}{{\cal D}}

\newcommand{\bc}{{\bf c}}
\newcommand{\ba}{{\bf a}}
\newcommand{\bb}{{\bf b}}

\newcommand{\bbf}{{\bf f}}

\newcommand{\bh}{{\bf h}}

\newcommand{\bp}{{\bf p}}
\newcommand{\bq}{{\bf q}}

\newcommand{\bs}{{\bf s}}
\newcommand{\bx}{{\bf x}}
\newcommand{\bu}{{\bf u}}
\newcommand{\bv}{{\bf v}}
\newcommand{\bw}{{\bf w}}
\newcommand{\bi}{{\bf i}}

\newcommand{\by}{{\bf y}}
\newcommand{\bA}{{\bf A}}

\newcommand{\bG}{{\bf G}}
\newcommand{\bJ}{{\bf J}}
\newcommand{\bK}{{\bf K}}
\newcommand{\bH}{{\bf H}}
\newcommand{\bL}{{\bf L}}
\newcommand{\bM}{{\bf M}}

\newcommand{\bI}{{\bf I}}

\newcommand{\bW}{{\bf W}}
\newcommand{\bU}{{\bf U}}
\newcommand{\bY}{{\bf Y}}

\newcommand{\bDelta}{{\utwi{\Delta}}}

\newcommand{\bGamma}{{\utwi{\Gamma}}}
\newcommand{\bLambda}{{\utwi{\Lambda}}}

\newcommand{\bone}{{\bf 1}}

\newcommand{\rev}[1]{{\color{black} #1}}
\newcommand{\revv}[1]{{\color{black} #1}}
\newcommand{\revvv}[1]{{\color{black} #1}}

\newcommand{\jay}{\jmath}
\newcommand{\conj}{\overline}
\newcommand{\reals}{\mathbb{R}}
\newcommand{\comps}{\mathbb{C}}

\newcommand{\sfT}{\textsf{T}}

\DeclarePairedDelimiterX{\norm}[1]{\lVert}{\rVert}{#1}

\usepackage{epstopdf}

\newtheorem{lemma}{Lemma}
\newtheorem{theorem}{Theorem}

\theoremstyle{definition}
\theoremstyle{definition}
\newtheorem{remark}{Remark}

\begin{document}

\captionsetup{font=small}

\IEEEoverridecommandlockouts

%------------------------------------------------------------------------------
% Title.
%------------------------------------------------------------------------------
%\title{Load-Flow Solutions %in Three-Phase Distribution Networks 
%with Mixed Delta-wye Connections: Existence, Uniqueness, and Linear Models} 

\title{Load-Flow in Multiphase Distribution Networks: Existence, Uniqueness, \rev{Non-Singularity} and Linear Models}

\author{Andrey Bernstein, %\emph{Member, IEEE}, 
Cong Wang,  %\emph{Member, IEEE}, 
Emiliano Dall'Anese,  %\emph{Member, IEEE}, 
Jean-Yves Le Boudec,  %\emph{Fellow, IEEE}, 
and Changhong Zhao%,  \emph{Member, IEEE}
\thanks{Andrey Bernstein, Emiliano Dall'Anese, and Changhong Zhao are with the National Renewable Energy Laboratory (NREL), Golden, CO, USA. %Emails: name.lastname@nrel.gov. 
Cong Wang and Jean-Yves Le Boudec are with \'{E}cole Polytechnique F\'{e}d\'{e}rale de Lausanne (EPFL), Lausanne, Switzerland. %Emails: name.lastname@epfl.ch.
}
\thanks{The work of A. Bernstein, E. Dall'Anese, and C. Zhao  was supported by the U.S. Department of Energy under Contract No. DE-AC36-08GO28308 with the National Renewable Energy Laboratory; funds provided by the Advanced Research Projects Agency-Energy (ARPA-E) under the Network Optimized Distributed Energy Systems (NODES) program.  The U.S. Government retains and the publisher, by accepting the article for publication, acknowledges that the U.S. Government retains a nonexclusive, paid-up, irrevocable, worldwide license to publish or reproduce the published form of this work, or allow others to do so, for U.S. Government purposes.

We would like to thank Fei Ding for her extensive help and support with the OpenDSS software.}
}

\maketitle

\begin{abstract}
This paper considers unbalanced multiphase distribution systems \rev{with generic topology and different load models,} and \rev{extends the $Z$-bus iterative} load-flow algorithm based on a fixed-point interpretation of the AC load-flow equations.  Explicit conditions for existence and uniqueness of load-flow solutions are presented. These conditions also guarantee  convergence of the load-flow algorithm  to the unique solution. The proposed methodology %broadens the well-established $Z$-bus iterative method, and it 
is applicable to generic systems featuring (i) wye connections; (ii) ungrounded delta connections; (iii) a combination of wye-connected and delta-connected sources/loads; and, (iv) a combination of line-to-line and line-to-grounded-neutral devices at the secondary of distribution transformers. \rev{Further, a sufficient condition for the non-singularity of the load-flow Jacobian is proposed. Finally,} linear load-flow models are derived, and their approximation accuracy is analyzed. 
% The development of approximate linear models is motivated by the need of computationally-affordable optimization and control applications -- from advanced distribution management systems settings to online optimal power flow routines.
Theoretical results are corroborated through experiments on IEEE test feeders.   
\end{abstract}

%\begin{IEEEkeywords}
%Power system operations, power system computational analysis, power-flow methods, linear power-flow models, multiphase distribution networks
%\end{IEEEkeywords} 

\section{Introduction} 
%\textcolor{blue}{Explain better when you have both line-line and line-neutral connections at the same node: 1) when multiple transformers are bundled together at an electrical node for model reduction purposes; 2) when we model also the secondary of the network and we have line-line and line-neutral loads. Do not use delta and wye terms; replace with line-line and line-neutral. Look: 342-Node Low Voltage Network Test System.}. 

%\textcolor{blue}{Explain better the advantages compared to Gatsis, when we have disjoint line-line and line-neutral; larger convergence region!}

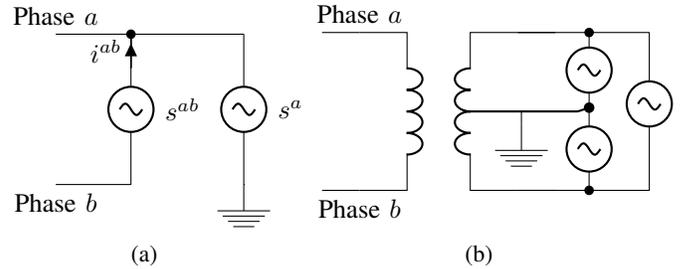
\begin{figure}[t]
\centering
\begin{subfigure}[b]{0.205\textwidth}
\begin{circuitikz}[american voltages,scale=1,transform shape]
\draw
  (0,0) node[above] {Phase $a$} to [short, -] (1,0)
  to [/tikz/circuitikz/bipoles/length=1cm, vsourcesin, l = $s^{ab}$] (1, -2)
  (1, -0.45) to [short, i=$i^{ab}$] (1, 0)
  (1, 0) to [short, *-] (2.5, 0)
  to [/tikz/circuitikz/bipoles/length=1cm, vsourcesin, l = $s^{a}$] (2.5, -2)
  to node[ground] {} (2.5, -2)
  (0,-2) node[below] {Phase $b$} to [short, -] (1,-2)
  ;
\end{circuitikz}
\caption{} \label{fig:deltawye1}
\end{subfigure}%
\hspace*{0.2cm}
\centering
\begin{subfigure}[b]{0.25\textwidth}
\begin{circuitikz}[american,scale=1,transform shape]
\draw (0,0) node [transformer](T){}  % reminded by @PaulGessler, thanks.
      (T.A1) node[above] {Phase $a$}
      (T.A2) node[below] {Phase $b$}
      (T.B1) node[above] {} 
      (T.B2) node[below] {}
      (T.base) node{};
\draw (T.A1) --++(-0.5,0);
\draw (T.A2) --++(-0.5,0);
\draw (T.B1) --++(0.5,0) to[open](4,0);  
\draw (T.B2) --++(0.5,0) to[open](3,-2);
\draw(T.A1) to[open](T.A2);
\draw(T.B2) to[open](T.B1);
% 2 new lines for neutral line on the secondary side.
\draw[thick] ($(T.B1)!0.5!(T.B2)-(0.7,0)$)--node[pos=0.5,above,inner sep=0pt](n){}++ (1.5,0) to [short, -*] (2, -1);
\draw  (n) -- ++ (0,-0.2)node[ground]{}
(T.B1) to [short, -*] (2, 0)
to [/tikz/circuitikz/bipoles/length=1cm, vsourcesin] (2, -1)
to [/tikz/circuitikz/bipoles/length=1cm, vsourcesin] (2, -2 |- T.B2)
to [short, *-] (T.B2)
(2, -2 |- T.B2) to [short, -] (2.8, -2 |- T.B2)
to [short, -] (2.8, -1.9)
to [/tikz/circuitikz/bipoles/length=1cm, vsourcesin] (2.8, 0)
to [short, -] (2, 0)
;
\end{circuitikz}
\caption{} \label{fig:deltawye2}
\end{subfigure}
\caption{Examples of a multiphase point of connection (for illustration simplicity, only two phases are shown): (a) Combination of wye-connected and delta-connected loads/sources at the primary of the feeder (due to e.g., network reduction procedures). (b) Combination of line-to-line loads/sources and  line-to-grounded-neutral devices at the secondary of distribution transformer. }
\label{fig:deltawye}
\vspace{-.5cm}
\end{figure}

Load-flow analysis
%which expresses the link between complex node voltages and complex nodal power injections, 
is a fundamental  task in power system theory and applications. 
In this paper, we consider a load-flow problem for a multiphase distribution network. The network has a generic topology (it can be either radial or meshed), it has a single slack bus with voltages that are fixed and known, and it features multiphase $PQ$ buses.  At each multiphase bus, the model of the  distribution system can have: (i) grounded wye-connected loads/sources; (ii) ungrounded delta connections; (iii) a combination of wye-connected and delta-connected loads/sources; or, (iv) a combination of line-to-line and line-to-grounded-neutral devices at the secondary of distribution transformers~\cite{Kerstingbook}. Models (i)--(iii) pertain to settings when the network model is limited to (aggregate) nodal power injections at the primary side of distribution transformers. Particularly, the combined  model (iii) can be utilized when different distribution transformers with either delta and/or wye primary connections are bundled together at one bus for network reduction purposes (e.g., when two transformers are connected through a short low-impedance line); see Figure \ref{fig:deltawye}(a) for an illustration. Load model (iv) is common in, e.g., North America for commercial buildings and residential customers, and it can be utilized when the network model includes the secondary of the distribution transformers\footnote{\rev{We note that models (iii) and (iv) are the same in terms of the mathematical formulation. However, from the practical point of view, model (iv) reflects an actual mode of connection on the secondary side of the distribution transformer, whereas model (iii) pertains to the case where different distribution transformers are lumped in the same bus for network reduction purposes.}}; see an illustrative example in Figure~\ref{fig:deltawye}(b) and low-voltage test feeders available in the literature (e.g., the IEEE 342-Node Low-Voltage Test System). Settings with only line-line or line-ground connections at the secondary are naturally subsumed by model (iv). 

%************* delta-wye figure
%***********************

Due to the nonlinearity of the AC load-flow equations, the existence and uniqueness of the solution to the load-flow problem is not guaranteed \rev{globally}. \rev{In fact, it is well known that the load-flow problem might have multiple solutions, as shown, e.g., in} \cite{multSol1,multSol2,multSol3}. Recently, solvability of 
lossless load-flow equations was investigated in~\cite{Simpson16}. Focusing on the exact AC load-flow equations, several efforts investigated explicit conditions for  existence  and uniqueness of the (high-voltage) solution \rev{within a given domain} in balanced distribution networks  \cite{Bolognani15, yu2015, congLF} as well as in the more realistic case of unbalanced three-phased networks \cite{congThreePhase, ZIP}. 

This paper examines the load-flow problem for multiphase distribution systems \rev{with any topology} and load models (i)--(iv), and outlines a load-flow iterative solution method that broadens the classical $Z$-bus  methodologies \cite{chen1991,Borzacchiello}. The iterative algorithm is obtained by leveraging the fixed-point interpretation of the nonlinear AC load-flow equations in \cite{congLF}. 
\rev{The specific formulation of the load-flow problem allows us to obtain explicit theoretical conditions that guarantee the existence of the load-flow solution that is unique in a domain that is analytically characterized}. Under these conditions, it is shown that the iterative algorithm achieves this unique solution. Compared to existing  methods and analysis, the contribution is \rev{threefold}:
\begin{itemize}

\item When only the load models (i)--(ii) are utilized (and for settings with only line-line or line-ground connections at the secondary), the analytical conditions for convergence presented in this paper improve upon existing methods \cite{congThreePhase, ZIP} by providing an enlarged set of power profiles that guarantee convergence. 

\item The methods and analysis outlined in \cite{congThreePhase, ZIP} are not applicable when the load models (iii) and (iv) are utilized. On the other hand, this paper provides a unified load-flow solution method for general load models at both the primary and secondary sides of the distribution transformer. 
To the best of our knowledge, the only existing freely available load-flow solver for networks with all models (i)--(iv) is part of the OpenDSS platform \cite{openDSS}. \rev{In fact, the algorithm utilized there \cite{openDSSsol} is based on a fixed-point iteration -- similar to our method, although not identical.  Our methodology can be conceivably extended to analyze the convergence properties of \cite{openDSSsol}.}

\item \rev{A sufficient condition for the non-singularity of the load-flow Jacobian is presented. Moreover, we show that the solutions guaranteed by our conditions satisfy the non-singularity of the load-flow Jacobian.}
\end{itemize}
\rev{We note the iterative solution method proposed in this paper is similar to the fixed-point MANA method in \cite{kocarMANA}, although no convergence results are provided in  \cite{kocarMANA}. The iterations in  \cite{kocarMANA} are not explicitly formulated in terms of voltage phasors, and hence it might be hard to analyze  its convergence properties using the tools outlined in this paper. It is also worth noticing that \cite{kocarMANA} does not consider  delta-connected loads.}

The paper then presents and analyzes two approximate load-flow models\footnote{\revvv{A follow up work appeared in the 7th IEEE International Conference on Innovative Smart Grid Technologies, ISGT Europe 2017, under the title ``Linear Power-Flow Models in Multiphase Distribution Networks.''  The follow up paper develops additional linear models for power flow at the substation and line currents, and outlines some applications.}} to relate voltages and complex power injections through an approximate linear relationship.
\rev{The first model is based on a standard application of the first-order Taylor (or tangent plane) local approximation, whereas the second model is directly based on our fixed-point formulation of the load-flow equations. The latter model provides a non-local approximation of the load-flow solution and is in the spirit of the previously proposed linear model for balanced networks in \cite{Bolognani15}.}

The development of approximate linear models is motivated by the need of computationally-affordable optimization and control applications -- from advanced distribution management systems settings to online and distributed optimization routines. For example, the nonlinearity of the (exact) AC load-flow equations poses significant difficulties in solving AC optimal power flow (OPF) problems \cite{Low-Convex,swaroop2015linear}. Typical approaches involve convex relaxation methods (e.g., semidefinite program \cite{Low-Convex}) or a linearization of the load-flow equations \cite{christ2013sens,sairaj2015linear,bolognani2015linear}. 
\rev{For multiphase unbalanced settings, linear load-flow models have been recently  proposed in \cite{garcesLinear,ahmadiLinear,kekatosLinear}. 
In particular, the method in \cite{garcesLinear} is based on the Taylor expansion of complex-valued functions; however, the extension to the general unbalanced case with a combination of delta and wye connections is not presented.  In \cite{ahmadiLinear},  a curve-fitting technique is used to fit a linear model to the non-linear load-flow equations. In order to treat the delta loads, they are translated into equivalent wye loads; therefore, the method cannot be used explicitly in the optimization settings where the power consumed/produced by the delta loads constitutes a control variable. In \cite{kekatosLinear}, an extension of the LinDistFlow model to a multiphase setting is proposed; however, the method is only applicable to radial grids, and no delta loads are considered. Moreover, no theoretical bounds on the approximation error are provided in \cite{garcesLinear,ahmadiLinear,kekatosLinear}.}

Approximate linear models have been recently utilized to develop  real-time OPF solvers for distribution systems \cite{commelec1,opfPursuit}. The  methodology proposed in the present paper is applicable to generic multiphase networks, and it thus can be utilized to broaden the applicability of~\cite{swaroop2015linear,commelec1,opfPursuit}.

\iffalse
\begin{circuitikz}[american]
\draw (0,0) node [transformer](T){}  % reminded by @PaulGessler, thanks.
      (T.A1) node[above] {A1}
      (T.A2) node[below] {A2}
      (T.B1) node[above] {B1} 
      (T.B2) node[below] {B2}
      (T.base) node{K};
\draw (T.A1) --++(-2,0);
\draw (T.A2) --++(-2,0);
\draw (T.B1) --++(2,0) to[D, v=${V_\gamma=0.7}$, i>_=](5,0);        
\draw (T.B2) --++(2,0) to[D, v=${V_\gamma=0.7}$ ,i>_=](5,-2.1);
\draw(T.A1) to[open,v<={$240V_{rms}$}](T.A2);
\draw(T.B2) to[open,v>=$$](T.B1);
% 2 new lines for neutral line on the secondary side.
\draw[thick] ($(T.B1)!0.515!(T.B2)-(0.7,0)$)--node[pos=0.5,above,inner sep=0pt](n){$12V_{rms}$}++ (3,0);
\draw  (n) -- ++ (0,-0.3)node[ground]{};
\end{circuitikz}
\fi

\begin{table}[t]
\caption{\rev{Nomenclature}}
\label{table:nom}
%\begin{center}
\begin{tabular}{ll}
\hline \\
$N$: &number of $PQ$ buses\\
$j \in \{1, \ldots, N\}$: &index of a $PQ$ bus;\\
$\bs^Y_j = (s^a_j, s^b_j, s^c_j)^\sfT$: &grounded wye sources at bus $j$;\\
$\bs^{\Delta}_j = (s^{ab}_j, s^{bc}_j, s^{ca}_j)^\sfT$: &delta sources at bus $j$;\\
$\bv_j = (v^a_j, v^b_j, v^c_j)^\sfT$: &phase-to-ground voltages at bus $j$;\\
$\bi_j = (i^a_j, i^b_j, i^c_j)^\sfT$: &phase net current injections at bus $j$;\\
$\bi^{\Delta}_j = (i^{ab}_j, i^{bc}_j, i^{ca}_j)^\sfT$: &phase-to-phase currents at bus $j$;\\ 
$\bv_0 = (v_0^a, v_0^b, v_0^c)^\sfT$: &voltages at the slack bus;\\
$\bv = (\bv_1^\sfT, \ldots, \bv_N^\sfT)^\sfT$: &voltages at $PQ$ buses;\\ 
$\bi = (\bi_1^\sfT, \ldots, \bi_N^\sfT)^\sfT$: &current injections at $PQ$ buses;\\ 
$\bi^{\Delta} = ((\bi_1^\Delta)^\sfT, \ldots, (\bi_N^\Delta)^\sfT)^\sfT$: &phase-to-phase currents at $PQ$ buses;  \\  
$\bs^Y = ((\bs^Y_1)^\sfT, \ldots, (\bs^Y_N)^\sfT)^\sfT$: &wye sources at  $PQ$ buses;  \\ 
$\bs^{\Delta} = ((\bs_1^\Delta)^\sfT, \ldots, (\bs_N^\Delta)^\sfT)^\sfT$: &delta sources at  $PQ$ buses; \\
$\bY$: &multiphase admittance matrix;\\
$\bY_{LL}$: &$\bY$ matrix with slack bus removed;\\
$\bH$: &transformation block-diagonal matrix \\
& (phase-ground $\rightarrow$ phase-phase);\\
$\bw$:  &zero-load voltage profile;\\
$\xi^Y(\bs), \xi^\Delta(\bs), \xi(\bs)$: &norms that are used to define regions\\
&of existence and uniqueness;\\
$\alpha(\bv), \beta(\bv), \gamma(\bv)$: &voltage quantities that are used to define\\
&regions of existence and uniqueness;\\
$\bp^Y, \bq^Y, \bp^\Delta, \bq^\Delta$: &active and reactive power injections;\\
$\bx^Y = ((\bp^Y)^\sfT, (\bq^Y)^\sfT)^\sfT$: &stacked vector of wye-injections;\\
$\bx^\Delta = ((\bp^\Delta)^\sfT, (\bq^\Delta)^\sfT)^\sfT$: &stacked vector of delta-injections. \\
\hline
\end{tabular}
%\end{center}
\end{table}

%The paper is organized as follows. In Section \ref{sec:problem}, we outline the load-flow problem and describe the proposed solution method. In Section \ref{sec:existence}, we give our main result on existence and uniqueness. In Section \ref{sec:linear}, we propose two different linear power-flow models. In Section \ref{sec:numerical}, we present the results of numerical evaluation of the proposed methods using standard IEEE benchmarks. Finally, we conclude in Section \ref{sec:conclusion}.

\rev{\section{Nomenclature and Notation} }
Upper-case (resp. lower-case) boldface letters are used for matrices (resp. column vectors); $(\cdot)^\sfT$ for transposition; %$(\cdot)^*$ complex-conjugate; and, $(\cdot)^\sfH$ complex-conjugate transposition; $\Re\{\cdot\}$ and $\Im\{\cdot\}$ denote the real and
  %imaginary parts of a complex number, respectively; $\mathrm{j} := \sqrt{-1}$ the imaginary unit; 
  $|\cdot|$ for the absolute value of a number or the component-wise absolute value of a vector or a matrix; and the %non-italic 
  letter $\jay$  for $\jay := \sqrt{-1}$. %Let $\cA \times \cB$ denote the Cartesian product of sets $\cA$ and $\cB$.  %For $x \in \mathbb{R}$, function $[x]_+$ is defined as $[x]_+ := \max\{0,x\}$. 
  For a complex number $c \in \comps$, $\Re\{c\}$ and $\Im\{c\}$ denote its real and imaginary part, respectively; and $\conj{c}$ denotes the conjugate of $c$. For an $N \times 1$ vector $\bx \in \mathbb{C}^N$,  $\|\bx\|_\infty := \max(|x_1|...|x_n|)$, $\|\bx\|_1 := \sum_{i = 1}^N|x_i|$, and $\diag(\bx)$ returns an $N \times N$ matrix with the elements of $\bx$ in its diagonal. %; and, $\bx \succeq \by$ implies that the inequality $x_i \geq y_i$ is enforced for all the vector entries $i = 1, \ldots, N$.  
  %The spectral radius $\rho(\cdot)$ is defined for an $N\times N$ matrix $\bA$ and corresponding eigenvalues $\lambda_1 ... \lambda_N$ as $\rho(\bA) := \max(|\lambda_1|, ... ,|\lambda_N|)$. 
  For an $M\times N$ matrix $\bA \in \comps^{M \times N}$, the $\ell_\infty$-induced norm is defined as $||\bA||_\infty = \max_{i = 1, \ldots, M} \sum_{j = 1}^N |(\bA)_{ij}|$. Finally, for a vector-valued map $\bx: \by\in\reals^{N \times 1}\to \bx(\by)\in \comps^{M \times 1}$, we let $\partial \bx / \partial \by$ denote %the Jacobian matrix; that is,  
  the $M \times N$ complex matrix with entries $(\partial \bx / \partial \by)_{ik} =  \partial x_i / \partial y_k = \partial \Re\{x_i\}/ \partial y_k + \jay \partial \Im\{x_i\}/ \partial y_k $, $i =  1,\ldots,M$, $k = 1,\ldots,N$.
  %and the spectral norm is defined as $||\bA||_2 := \sqrt{\lambda_{max}(\bA^*\bA)}$, where $\lambda_{max}$ denotes maximum eigenvalue. Finally, $\bI_N$ denotes the $N \times N$ identity matrix.
\rev{
Nomenclature is given in Table \ref{table:nom}. Where possible, the definitions are also recalled upon use in the text.}

\section{Problem Formulation} \label{sec:problem}

For notational simplicity, the framework is outlined for three-phase systems; we describe in Remark \ref{rem:missingPhases} below how to apply the analysis to the general multiphase case (as we do in the numerical examples in Section \ref{sec:feed123}). Consider a generic three-phase distribution network with one slack bus and $N$ three-phase $PQ$ buses. With reference to the illustrative example in Figure \ref{fig:deltawye}, let $\bs^Y_j := (s^a_j, s^b_j, s^c_j)^\sfT$ denote the vector of grounded wye sources at bus $j$, where $s^\phi_j \in \comps$ denotes the net complex power injected on phase $\phi$. Similarly, let $\bs^{\Delta}_j := (s^{ab}_j, s^{bc}_j, s^{ca}_j)^\sfT$ denote the power injections of  delta-connected sources. With a slight abuse of notation, $\bs^Y_j$ and $\bs^{\Delta}_j$ will represent \revvv{line-ground and line-line} connections, respectively, when bus $j$ corresponds to the secondary side of the distribution transformer (this notational choice allows us not to introduce additional symbols). 

At bus $j$, the following set of equations relates voltages, currents, and powers:
%\[
%\left\lbrace
%\begin{array}{l}
%s^{ab}_j = (v^a_j - v^b_j) \conj{i^{ab}_j}, \\
%s^{bc}_j = (v^b_j - v^c_j) \conj{i^{bc}_j}, \\
%s^{ca}_j = (v^c_j - v^a_j) \conj{i^{ca}_j}, \\
%v^a_j (\conj{i^{ab}_j} - \conj{i^{ca}_j}) + s^a_j = v^a_j \conj{i^a_j}, \\
%v^b_j (\conj{i^{bc}_j} - \conj{i^{ab}_j}) + s^b_j = v^b_j \conj{i^b_j}, \\
%v^c_j (\conj{i^{ca}_j} - \conj{i^{bc}_j}) + s^c_j = v^c_j \conj{i^c_j},
%\end{array}
%\right.
%\]
\begin{align*}
&s^{ab}_j = (v^a_j - v^b_j) \conj{i^{ab}_j}, &v^a_j (\conj{i^{ab}_j} - \conj{i^{ca}_j}) + s^a_j = v^a_j \conj{i^a_j},\\
&s^{bc}_j = (v^b_j - v^c_j) \conj{i^{bc}_j}, &v^b_j (\conj{i^{bc}_j} - \conj{i^{ab}_j}) + s^b_j = v^b_j \conj{i^b_j},\\
&s^{ca}_j = (v^c_j - v^a_j) \conj{i^{ca}_j}, &v^c_j (\conj{i^{ca}_j} - \conj{i^{bc}_j}) + s^c_j = v^c_j \conj{i^c_j},
\end{align*}
where $\bv_j = (v^a_j, v^b_j, v^c_j)^\sfT$, $\bi_j = (i^a_j, i^b_j, i^c_j)^\sfT$, and $\bi^{\Delta}_j = (i^{ab}_j, i^{bc}_j, i^{ca}_j)^\sfT$ collect the phase-to-ground voltages $\{v^\phi_j\}_{\phi \in \{a,b,c\}}$, phase net current injections $\{i^\phi_j\}_{\phi \in \{a,b,c\}}$, and phase-to-phase currents $\{i^{\phi \phi^\prime}_j\}_{\phi, \phi^\prime \in \{a,b,c\}}$ (for delta connections and line-line connections) of node $j$, respectively.

We next express the set of load-flow equations in vector-matrix form. To this end, let $\bv_0 := (v_0^a, v_0^b, v_0^c)^\sfT$ %and $\bi_0 := (i_0^a, i_0^b, i_0^c)^\sfT$
%, and $\bs_0 := (s_0^a, s_0^b, s_0^c)^\sfT$ 
denote the complex vector collecting the three-phase voltages % and currents
%, and power injection 
at the slack bus (i.e., the substation). Also, let $\bv := (\bv_1^\sfT, \ldots, \bv_N^\sfT)^\sfT$, $\bi := (\bi_1^\sfT, \ldots, \bi_N^\sfT)^\sfT$, $\bi^{\Delta} := ((\bi_1^\Delta)^\sfT, \ldots, (\bi_N^\Delta)^\sfT)^\sfT$, $\bs^Y := ((\bs^Y_1)^\sfT, \ldots, (\bs^Y_N)^\sfT)^\sfT$, and $\bs^{\Delta} := ((\bs_1^\Delta)^\sfT, \ldots, (\bs_N^\Delta)^\sfT)^\sfT$ be the vectors in $\comps^{3N}$ collecting the respective electrical quantities of the $PQ$ buses. The load-flow problem is then defined as solving for $\bv$ (and $\bi^{\Delta}$) in the following set of equations, where $\bs^Y$, $\bs^{\Delta}$, and $\bv_0$ are given:
\begin{subequations} \label{eqn:lf}
\begin{align} 
\diag\left(\bH^\sfT \conj{\bi^{\Delta}}\right)\bv + \bs^Y = \diag(\bv) \conj{\bi}, \label{eqn:lf_balance}\\
\bs^{\Delta} = \diag\left(\bH \bv \right)\conj{\bi^{\Delta}}, \label{eqn:lf_delta}\\
\bi = \bY_{L0} \bv_0 + \bY_{LL} \bv. \label{eqn:lf_i} %\\
%\bs_0 =  \diag(\bv_0)\left (\conj{\bY}_{00} \conj{\bv}_0 + \conj{\bY}_{0L} \conj{\bv} \right). \label{eqn:s0} 
\end{align}
\end{subequations}
In~\eqref{eqn:lf}, $\bY_{00} \in \comps^{3 \times 3}, \bY_{L0} \in \comps^{3N \times 3}, \bY_{0L} \in \comps^{3 \times 3N}$, and $\bY_{LL} \in  \comps^{3N \times 3N}$ are the submatrices of the three-phase admittance matrix
\begin{equation}
\bY := 
\begin{bmatrix}
\bY_{00} & \bY_{0L} \\
\bY_{L0} & \bY_{LL}
\end{bmatrix} \in \comps^{3(N+1) \times 3(N+1)},
\end{equation}
which can be formed from the topology of the network, the $\pi$-model of the transmission lines, and other passive network devices, as shown in, e.g.,~\cite{Kerstingbook}; and $\bH$ is a $3N \times 3N$ block-diagonal matrix defined by
\begin{equation}
\bH :=
\begin{bmatrix}
\bGamma & & \\
& \ddots & \\
& & \bGamma
\end{bmatrix}
, \quad \bGamma := 
\begin{bmatrix}
1 & -1 & 0\\
0 & 1 & -1\\
-1 & 0 & 1
\end{bmatrix}.
\end{equation}
\rev{
In more detail, \eqref{eqn:lf_balance} follows from the Kirchoff's current law at the buses, \eqref{eqn:lf_delta} relates power injections and currents for the delta-connected loads/sources, and \eqref{eqn:lf_i} relates nodal current injections and voltages through Ohm's law.
}

By simple algebraic manipulations, \rev{$\bi^{\Delta} $ can be eliminated from the set \eqref{eqn:lf}, and  }the solution $\bv$ can be found from the following fixed-point equation:
\begin{equation} \label{eqn:fp}
\begin{split}
\bv & = \mathbf{G}_{\bs^{Y}\bs^\Delta}(\bv) \\ & := \bw + \bY_{LL}^{-1}\left(\diag(\conj{\bv})^{-1}\conj{\bs}^{Y} + \bH^\sfT \diag\left(\bH \conj{\bv} \right)^{-1} \conj{\bs}^{\Delta} \right),
\end{split}
\end{equation}
where\footnote{It was shown in \cite{congThreePhase,sairaj2015linear} that $\bY_{LL}$ is invertible for most practical cases of three-phase distribution networks.} $\bw := -\bY_{LL}^{-1}\bY_{L0}\bv_0$ is the zero-load voltage.

\rev{We note that the benefit of the proposed load-flow formulation \eqref{eqn:fp} is that it can be analyzed theoretically using the Banach fixed-point theory, as presented in the next section.}

\revv{Before proceeding, we recall the notion of \emph{non-singularity} associated with a load-flow solution. Note that \eqref{eqn:lf} defines an explicit mapping from the state vector $(\bv^\sfT, (\bi^\Delta)^\sfT)^\sfT \in \comps^{6N}$ to the vector of power injections $\bs := ((\bs^Y)^\sfT, (\bs^\Delta)^\sfT)^\sfT \in \comps^{6N}$.  Let $\bx := \left ( \Re\{\bs^Y\}^\sfT, \Im\{\bs^Y\}^\sfT, \Re\{\bs^\Delta\}^\sfT, \Im\{\bs^\Delta\}^\sfT \right )^\sfT$ denote the \emph{real-valued} vector that collects the active and reactive power injections of wye and delta sources. 
Similarly, let $\by := \left(\Re\{\bv\}^\sfT, \Im\{\bv\}^\sfT, \Re\{\bi^\Delta\}^\sfT, \Im\{\bi^\Delta\}^\sfT \right)^\sfT$ denote the \emph{real-valued} vector of the state variables. Then, the load-flow equations can be written as
\begin{equation} \label{eqn:real_lf}
\bx = \bh (\by),
\end{equation}
where $\bh: \reals^{12N} \rightarrow  \reals^{12N}$ is the mapping defined explicitly by \eqref{eqn:lf}. Let $\bJ(\by)$ be the Jacobian matrix of this mapping, i.e.,
%
%To this end, let $\bJ$ be the Jacobian matrix of the load-flow equation  \eqref{eqn:lf}, namely
$
(\bJ)_{ij} = \frac{ \partial (\bh)_i}{ \partial (\by)_j},~i, j\in\{ 1, \ldots, 12N\}$. We say that a given state vector $\by$ (and hence, the corresponding complex-valued vector $(\bv^\sfT, (\bi^\Delta)^\sfT)^\sfT$) is non-singular if the Jacobian matrix $\bJ(\by)$ is invertible. A pair $(\bv, \bs)$ is non-singular if the corresponding state vector $\left(\bv, \bi^\Delta := \diag^{-1}\left(\bH \conj{\bv} \right)\conj{\bs}^{\Delta}\right)$ is non-singular. The non-singularity property represents a sufficient condition for the (static) \emph{voltage stability} of the operating point (see, e.g., \cite{Kerstingbook}).
}

\begin{remark} \label{rem:missingPhases} %For notational simplicity, the proposed framework is outlined for three-phase systems. 
Observe that \eqref{eqn:fp} can be straightforwardly utilized in cases when a network features a mix of three-phase, two-phase, and single-phase buses. In particular, in that case, the vectors $\bv$, $\bs^Y$, and $\bw$ collect their corresponding electrical quantities only for existing phases; the vector $\bs^\Delta$ collects the existing phase-to-phase injections; and the matrix $\bH$ contains rows that correspond to the existing phase-to-phase connections. For example, if a certain bus has only a single $ab$ connection, it will only contain a row with $(1, -1, 0)$  for that bus. To be more precise, $\bH$ is $N^\Delta \times N^{phases}$ matrix, where $N^\Delta$ is the total number of phase-to-phase connections, and $N^{phases}$ is the total number of phases in all the buses. \rev{In the cases where there is no phase-to-phase connection in the network, the fixed-point formulation \eqref{eqn:fp} still holds after removing the term that involves $\bH,\bs^\Delta$.}
\end{remark}

\begin{remark}\label{rem:ZIPmodel}
For exposition simplicity, the  proposed method is outlined for the case of a constant-power load model. This is also motivated by recent optimization and control frameworks for distribution systems, where distributed energy resources as well as noncontrollable assets are (approximately) modeled as constant-PQ units~\cite{Low-Convex,swaroop2015linear,commelec1,opfPursuit,bolognani2015linear}. 
\rev{The extension of the results in the present paper to a more general ZIP load model is possible using the methodology of \cite{ZIP}; however, it is out of the scope of this paper.}
\end{remark}

\section{Existence, Uniqueness, \revv{and Non-Singularity}} \label{sec:existence}

The fixed-point equation~\eqref{eqn:fp} leads to an iterative procedure wherein the vector of voltages is updated as:   
\begin{equation} \label{eqn:fp_iter}
\mathbf{v}^{(k+1)}=\mathbf{G}_{\mathbf{s}^\mathbf{Y}\mathbf{s}^\Delta}(\mathbf{v}^{(k)})
\end{equation}
with $\mathbf{v}^{(0)}$ a given initialization point, $k$  the iteration index, and $\mathbf{G}_{\mathbf{s}^\mathbf{Y}\mathbf{s}^\Delta}(\cdot)$ defined in~\eqref{eqn:fp}. \rev{In fact, iteration \eqref{eqn:fp_iter} can be viewed as an extension of the classic $Z$-bus method to the general setting considered in this paper.} Convergence of the iterative method~\eqref{eqn:fp_iter} is analyzed next.  

To this end, let $\bW := \diag (\bw)$, and $\bL := |\bH|$ be the component-wise absolute value of the matrix $\bH$.
\iffalse
Also, let
\begin{equation} \label{eqn:L}
\bL := 
\begin{bmatrix}
\bLambda & & \\
& \ddots & \\
& & \bLambda
\end{bmatrix}
, \quad \bLambda := 
\begin{bmatrix}
1 & 1 & 0\\
0 & 1 & 1\\
1 & 0 & 1
\end{bmatrix},
\end{equation}
\fi
Also, for $\bs := ((\bs^Y)^\sfT, (\bs^\Delta)^\sfT)^\sfT \in \comps^{6N}$ define
\begin{subequations} \label{eqn:xi}
\begin{align}
\xi^Y(\bs) & := \left\|\bW^{-1} \bY_{LL}^{-1} \bW^{-1} \diag(\bs^Y) \right\|_{\infty}, \\
\xi^\Delta(\bs) & := \left\|\bW^{-1} \bY_{LL}^{-1} \bH^\sfT \diag (\bL|\bw|)^{-1} \diag(\bs^\Delta) \right\|_{\infty}, \\
\xi(\bs) & := \xi^Y(\bs) + \xi^\Delta(\bs),
\end{align}
\end{subequations}
% \begin{align} \label{eqn:xi}
% \xi(\bs) := &\left\|\bW^{-1} \bY_{LL}^{-1} \conj{\bW}^{-1} \diag(\conj{\bs}^Y) \right\|_{\infty} \nonumber\\
% &\, + \left\|\bW^{-1} \bY_{LL}^{-1} \bH^\sfT \diag (\bL|\bw|)^{-1} \diag(\conj{\bs}^\Delta) \right\|_{\infty},
% \end{align}
where $|\bw|$ is the component-wise absolute value of the vector $\bw$, and $\|\bA\|_{\infty}$ is the induced $\ell_\infty$-norm of a complex matrix $\bA$.
%Note that $\xi(\bs)$ defines a norm on $\comps^{6N}$. 

%%%% old method
\iffalse
and set
\begin{equation} \label{eqn:alpha}
\alpha := \|\bL |\bw| \|_{\infty},
\end{equation}
where $|\bw|$ is the component-wise absolute value of the vector $\bw$. Let $\bs := ((\bs^Y)^\sfT, (\bs^\Delta)^\sfT)^\sfT \in \comps^{6N}$ and define
\begin{align} \label{eqn:xi}
\xi(\bs) := &\left\|\bW^{-1} \bY_{LL}^{-1} \conj{\bW}^{-1} \diag(\conj{\bs}^Y) \right\|_{\infty} \nonumber\\
&\, + \frac{1}{\alpha}\left\|\bW^{-1} \bY_{LL}^{-1} \bH^\sfT \diag(\conj{\bs}^\Delta) \right\|_{\infty},
\end{align}
where  for any complex matrix $\bA$, $\|\bA\|_{\infty}$ is the induced $\ell_\infty$ norm.
\fi
%%%%

\begin{lemma} \label{lem:norm}
$\xi(\bs)$ is a norm on $\comps^{6N}$.
\end{lemma}

The proof of Lemma~\ref{lem:norm} as well as other technical results are deferred to the Appendix. Finally, let 
%%%% old method
%\begin{equation} \label{eqn:rho}
%\rho(\bv) := \min \left \{  \min_j \frac{|(\bv)_j|} {|(\bw)_j|}, \, \frac{1}{\alpha}\min_j |(\bH \bv)_j| \right\}
%\end{equation}
\begin{subequations} \label{eqn:gamma}
\begin{align}
\alpha(\bv) & := \min_j \frac{|(\bv)_j|} {|(\bw)_j|} \\
\beta(\bv) & := \min_j \frac{|(\bH \bv)_j|}{(\bL|\bw|)_j} \\
\gamma(\bv) & := \min \left \{  \alpha(\bv), \, \beta(\bv) \right\} \label{eqn:gammaV}
\end{align}
\end{subequations}
% \begin{equation} \label{eqn:rho}
% \rho(\bv) := \min \left \{  \min_j \frac{|(\bv)_j|} {|(\bw)_j|}, \, \min_j \frac{|(\bH \bv)_j|}{(\bL|\bw|)_j} \right\}
% \end{equation}
We next present our main result on the solution of the fixed-point equation defined by \eqref{eqn:fp}.
\begin{theorem} \label{thm:exist_general}
Let $\hat{\mathbf{v}}$ be a given solution to the load-flow equations for a vector of power injections $\hat{\mathbf{s}}$. Consider some other candidate vector of power injections $\mathbf{s}$, and assume that there exists a $\rho\in(0,\gamma(\hat{\bv}))$, such that
\begin{equation} \label{eqn:cond1_gen}
\frac{\xi^Y(\mathbf{s}-\hat{\mathbf{s}})+\displaystyle\frac{\xi^Y(\hat{\mathbf{s}})}{\alpha(\hat{\bv})}\rho}{\alpha(\hat{\bv})-\rho}+\frac{\xi^\Delta(\mathbf{s}-\hat{\mathbf{s}})+\displaystyle\frac{\xi^\Delta(\hat{\mathbf{s}})}{\beta(\hat{\bv})}\rho}{\beta(\hat{\bv})-\rho}\leq\rho
\end{equation}
and
\begin{equation} \label{eqn:cond2_gen}
\frac{\xi^Y(\mathbf{s})}{(\alpha(\hat{\bv})-\rho)^2}+\frac{\xi^\Delta(\mathbf{s})}{(\beta(\hat{\bv})-\rho)^2}<1.
\end{equation}
Then, there exists a unique solution $\mathbf{v}$ in
\begin{equation} \label{eqn:D}
\cD_\rho(\hat{\bv}) := \{\bv: \, |(\bv)_j - (\hat{\bv})_j| \leq \rho |(\bw)_j|, \, j = 1\ldots3N \}
\end{equation}
to the load-flow equations with power injection $\mathbf{s}$. Moreover, this solution can be reached by iteration~\eqref{eqn:fp_iter}
%\begin{equation} \label{eqn:fp_iter}
%\mathbf{v}^{(k+1)}=\mathbf{G}_{\mathbf{s}^\mathbf{Y}\mathbf{s}^\Delta}(\mathbf{v}^{(k)})
%\end{equation}
initialized anywhere in $\mathcal{D}_\rho(\hat{\mathbf{v}})$.
\end{theorem}

The conditions of Theorem~\ref{thm:exist_general} may be computationally intensive as they require a  parameter scanning to find a proper value for $\rho$. In the following, we sacrifice the tightness of the inequalities \eqref{eqn:cond1_gen} and \eqref{eqn:cond2_gen} to obtain the following \emph{explicit conditions}.
\begin{theorem} \label{thm:exist_lf}
Let $\hat{\mathbf{v}}$ be a given solution to the load-flow equations with power injection $\hat{\mathbf{s}}$ satisfying:
\begin{equation} \label{eqn:cond1}
\xi(\hat{\bs}) < (\gamma(\hat{\bv}))^2, 
\end{equation}
\rev{where $\xi(\cdot)$ and  $\gamma(\cdot)$ are given in \eqref{eqn:xi} and \eqref{eqn:gamma}, respectively.}
Consider some other candidate power injections vector $\mathbf{s}$, and assume that
\begin{equation} \label{eqn:cond2}
\xi(\bs - \hat{\bs}) < \frac{1}{4} \left(\frac{(\gamma(\hat{\bv}))^2 - \xi(\hat{\bs})}{\gamma(\hat{\bv})} \right)^2.
\end{equation} %\vspace*{-1cm}
Let
\begin{subequations} \label{eqn:rhos}
\begin{align}
\rho^\ddagger(\hat{\bv},\hat{\bs}) & := \frac{1}{2} \left(\frac{(\gamma(\hat{\bv}))^2 - \xi(\hat{\bs})}{\gamma(\hat{\bv})} \right) \\
\rho^\dagger(\hat{\bv},\hat{\bs}, \bs) & := \rho^\ddagger(\hat{\bv},\hat{\bs}) - \sqrt{\left(\rho^\ddagger(\hat{\bv},\hat{\bs})\right)^2-\xi(\bs - \hat{\bs})}
\end{align}
\end{subequations}
Then:
\begin{enumerate}[(i)]
\item \revv{The operating point $(\hat{\bv}, \hat{\bs})$ is non-singular.}
\item There exists a unique load-flow solution $\bv$ in $\cD_\rho(\hat{\bv})$ defined in \eqref{eqn:D} with $\rho=\rho^\ddagger(\hat{\bv},\hat{\bs})$; %with complex power injection $\bs$;

\item This solution can be reached by iteration \eqref{eqn:fp_iter} starting from anywhere in $\cD_\rho(\hat{\bv})$ with $\rho=\rho^\ddagger(\hat{\bv},\hat{\bs})$;
\item The solution is located in $\cD_\rho(\hat{\bv})$ with $\rho=\rho^\dagger(\hat{\bv},\hat{\bs}, \bs)$;
\item The pair $(\bv, \bs)$ satisfies $\xi(\bs) < (\gamma(\bv))^2$; \revv{hence, it is non-singular}.
\end{enumerate}
\end{theorem}

\rev{Some comments about the above results follow:}
\begin{enumerate}[(a)]
\item If a solution to the load-flow problem $(\hat{\bv},\hat{\bs})$ is not always available, one can simply set  $\hat{\bv} = \bw$ and $\hat{\bs} = 0$ (with $\bw$ the zero-load voltage profile); see, e.g.,~\cite{congThreePhase, ZIP}. \rev{In such a case, condition \eqref{eqn:cond1} is trivially satisfied, and the existence and uniqueness is determined  based on \eqref{eqn:cond2}.
With respect to \cite{ZIP}, the main innovation is in the fact that our methodology allows to provide better conditions whenever a known load-flow solution is available. This setting is of particular practical interest in real-time control of  power networks, whereby a measurement of the state is available at every time step, and thus conditions can be refined to reflect the uniqueness in a domain around a given operating point. This property is absent in \cite{ZIP}, and consequently it is easy to find a situation where the conditions of the present paper are applicable, whereas the conditions of \cite{ZIP} are not; see Section \ref{sec:numerical} for examples.
\item Theorem \ref{thm:exist_lf} provides explicit sufficient conditions under which  conditions \eqref{eqn:cond1_gen} and \eqref{eqn:cond2_gen} of Theorem \ref{thm:exist_general} are satisfied.
Moreover, the particular conditions' formulation of Theorem \ref{thm:exist_lf} allows for a better \emph{localization} of the unique solution. Indeed, note that Theorem \ref{thm:exist_lf} provides two balls around a given load-flow solution. The first, bigger ball given by Theorem \ref{thm:exist_lf} (i) specifies the region of uniqueness in the voltage space; whereas the second, smaller ball given by Theorem \ref{thm:exist_lf} (iii) localizes this solution. An illustration is provided in Section \ref{sec:simpleillu}.
\item The explicit conditions \eqref{eqn:cond1} and \eqref{eqn:cond2} are useful in the OPF settings. More precisely, \eqref{eqn:cond2} can be utilized as explicit \emph{convex} constraint that ensures existence and uniqueness of the load-flow solution.
\item \revv{Part (v) of Theorem \ref{thm:exist_lf} suggests a successive application of our results, producing a sequence of non-singular load-flow solutions.}
%In the next section, we will show that $\xi(\bs) < (\gamma(\bv))^2$ is sufficient for the load-flow Jacobian to be non-singular.
\item The general multiphase networks can be treated using the method described in Remark \ref{rem:missingPhases}. For networks where there is no phase-to-phase connection, the correctness of the proposed theory is preserved by eliminating all terms and variables that involve $\bH,\bL,\bs^\Delta$. More precisely in those cases, we have $\xi(\bs)=\xi^Y(\bs)$ and $\gamma(\bv)=\alpha(\bv)$ in \eqref{eqn:xi},\eqref{eqn:gamma},\eqref{eqn:cond1},\eqref{eqn:cond2}, and \eqref{eqn:rhos}. In addition, we remove the second term on the left-hand side of \eqref{eqn:cond1_gen} and \eqref{eqn:cond2_gen}.
}
\end{enumerate}

\section{\revv{Linear Models}} \label{sec:linear}
In this section, we develop two methods to obtain  approximate representations of the AC load-flow equations \eqref{eqn:lf}, wherein the net injected powers and voltages are related through an approximate linear relationship. The first method is based on the first-order Taylor (FOT) expansion of the load-flow solution around a given point. FOT is therefore the best local linear approximator. 
%\rev{In this context, we also show that \eqref{eqn:cond1} is a sufficient condition for the non-singularity of the load-flow Jacobian.}
The second method is based on a single iteration of the fixed-point iteration \eqref{eqn:fp_iter} and it is hereafter referred to as fixed-point linearization (FPL). 

Let $\bp^Y := \Re\{\bs^Y\}$, $\bq^Y := \Im\{\bs^Y\}$, $\bp^\Delta := \Re\{\bs^\Delta\}$, $\bq^\Delta := \Im\{\bs^\Delta\}$, $\bx^Y := ((\bp^Y)^\sfT, (\bq^Y)^\sfT)^\sfT$,  and $\bx^\Delta := ((\bp^\Delta)^\sfT, (\bq^\Delta)^\sfT)^\sfT$ collect the active and reactive power injections. Also, let $|\bv|$ collect the voltage magnitudes. Our goal is to derive linear approximations to \eqref{eqn:lf} in the form 
\begin{subequations} \label{eqn:lin_lf}
\begin{align} 
\tilde{\bv} = \bM^Y \bx^Y + \bM^{\Delta} \bx^{\Delta} + \ba, \label{eqn:lin_v}\\
|\tilde{\bv}| = \bK^Y \bx^Y + \bK^{\Delta} \bx^{\Delta} + \bb, \label{eqn:lin_rho}
%\tilde{\bs}_0 =  \bG^Y \bx^Y + \bG^{\Delta} \bx^{\Delta} + \bc. \label{eqn:lin_s0}
\end{align}
for some matrices $\bM^Y, \bM^{\Delta} \in \comps^{3N \times 6N}$, $\bK^Y, \bK^{\Delta} \in \reals^{3N \times 6N}$, and vectors $\ba \in \comps^{3N}, \bb \in \reals^{3N}$.
\end{subequations}

\subsection{First-Order Taylor (FOT) Method}
\revv{Let $(\hat{\bv}, \hat{\bi}^\Delta,  \hat{\bs}^Y, \hat{\bs}^\Delta)$ be a given operating point satisfying \eqref{eqn:lf}, and let $\hat{\by}$ and $\hat{\bx}$ be the corresponding real-valued vectors. 
To obtain \eqref{eqn:lin_v},
we plug \eqref{eqn:lf_i} into \eqref{eqn:lf_balance}, and take partial derivatives of \eqref{eqn:lf_balance} and \eqref{eqn:lf_delta} with respect to $\bx^Y$ and $\bx^\Delta$:}
\begin{subequations} \label{eqn:part_der}
\begin{align} 
&\diag\left(\bH^\sfT \conj{\bi^{\Delta}}\right) \frac{\partial \bv}{\partial \bx^Y} + \diag (\bv) \bH^\sfT \frac{\partial \conj{\bi^{\Delta}}}{\partial \bx^Y}  + \bU  \notag\\
& \, = \diag(\bv) \conj{\bY}_{LL} \frac{\partial \conj{\bv}}{\partial \bx^Y} + \diag (\conj{\bY}_{L0} \conj{\bv}_0 + \conj{\bY}_{LL} \conj{\bv}) \frac{\partial \bv}{\partial \bx^Y},
 \\
&{\bf 0} = \diag\left(\bH \bv \right)\frac{\partial \conj{\bi^{\Delta}}}{\partial \bx^Y} + \diag(\conj{\bi^{\Delta}})\bH \frac{\partial \bv}{\partial \bx^Y}  , \\
&\diag\left(\bH^\sfT \conj{\bi^{\Delta}}\right) \frac{\partial \bv}{\partial \bx^\Delta} + \diag (\bv) \bH^\sfT \frac{\partial \conj{\bi^{\Delta}}}{\partial \bx^\Delta}  \notag\\
& = \diag(\bv) \conj{\bY}_{LL} \frac{\partial \conj{\bv}}{\partial \bx^\Delta} + \diag (\conj{\bY}_{L0} \conj{\bv}_0 + \conj{\bY}_{LL} \conj{\bv}) \frac{\partial \bv}{\partial \bx^\Delta},
 \\
&\bU = \diag\left(\bH \bv \right)\frac{\partial \conj{\bi^{\Delta}}}{\partial \bx^\Delta} + \diag(\conj{\bi^{\Delta}})\bH \frac{\partial \bv}{\partial \bx^\Delta} ,
\end{align}
\end{subequations}
where $\bU := (\bI_{3N}, \jay \bI_{3N}) \in \comps^{3N \times 6N}$ and $\bI_{3N} \in \reals^{3N \times 3N}$ is the identity matrix.  In this set of equations, set $\bv = \hat{\bv}$ and $\conj{\bi^{\Delta}} =  \diag(\bH \hat{\bv})^{-1}\hat{\bs}^\Delta$; the unknowns are the matrices $\frac{\partial \bv}{\partial \bx^Y}, \frac{\partial \conj{\bi^{\Delta}}}{\partial \bx^Y}, \frac{\partial \bv}{\partial \bx^\Delta},  \frac{\partial \conj{\bi^{\Delta}}}{\partial \bx^\Delta} \in \comps^{3N \times 6N}$. \revv{Model \eqref{eqn:lin_v} is then obtained by solving \eqref{eqn:part_der} and setting}
%\vspace*{-0.05cm}
\begin{align*}
\bM^Y  := \frac{\partial \bv}{\partial \bx^Y}, \quad
\bM^\Delta := \frac{\partial \bv}{\partial \bx^\Delta},
\end{align*}
%\vspace*{-0.385cm}
and
$
\ba := \hat{\bv} - \bM^Y \hat{\bx}^Y - \bM^\Delta \hat{\bx}^\Delta.
$

Observe that, in rectangular coordinates,  \eqref{eqn:part_der} is a set of \emph{linear} equations with the same number, $(12N)^2$, of real-valued equations and variables. 
\revv{In fact, \eqref{eqn:part_der} can be written as $\bJ (\hat{\by}) \frac{\partial \by}{\partial \bx} = \bI_{12N}$, where $\bJ(\cdot)$ is the Jacobian of the load-flow mapping  $\bh(\cdot)$ defined in \eqref{eqn:real_lf}, and $\bI_{12N} \in \reals^{12N \times 12N}$ is the identity matrix. Clearly, this equation has a unique solution if and only if $\bJ (\hat{\by})$ is invertible, namely $\hat{\by}$ is non-singular. Note that a sufficient condition for that is given by condition \eqref{eqn:cond1} of Theorem \ref{thm:exist_lf} (cf.~item (i) of that Theorem).}

\iffalse
We next give conditions under which the system of equations~\eqref{eqn:part_der} has a unique solution. The load-flow equations \eqref{eqn:lf} define an explicit mapping \rev{ $\bx = \bbf(\by)$  with
$\bx := \left((\bx^Y)^\sfT, (\bx^\Delta)^\sfT\right)^\sfT$ and $\by := \left(\Re\{\bv\}^\sfT, \Im\{\bv\}^\sfT, \Re\{\bi^\Delta\}^\sfT, \Im\{\bi^\Delta\}^\sfT \right)^\sfT$. 
Let $\bJ$ be the Jacobian matrix of this mapping, i.e.,
%
%To this end, let $\bJ$ be the Jacobian matrix of the load-flow equation  \eqref{eqn:lf}, namely
$
(\bJ)_{ij} = \frac{ \partial (\bbf)_i}{ \partial (\by)_j},~i, j\in\{ 1, \ldots, 12N\}$.
 Then \eqref{eqn:part_der} can be compactly written as $\bJ \frac{\partial \by}{\partial \bx} = \bI_{12N}$ with $\bI_{12N} \in \reals^{12N \times 12N}$ being the identity matrix.}

\begin{theorem} \label{thm:taylor}
Let $(\hat{\bv}, \hat{\bs})$ be a given operating point. 
\rev{
\begin{enumerate}[(i)]
\item The system \eqref{eqn:part_der} has a unique solution if and only if the Jacobian $\bJ$ is non-singular.
\item Under condition \eqref{eqn:cond1}, the Jacobian matrix $\bJ$ is non-singular and the system \eqref{eqn:part_der} has a unique solution.
\end{enumerate}
}
\end{theorem}

\rev{
Combining  (iv) of Theorem \ref{thm:exist_lf} and  (ii) of Theorem \ref{thm:taylor}, one has that the solution guaranteed in Theorem \ref{thm:exist_lf} automatically satisfies the non-singularity of the load-flow Jacobian.
}
\fi

To obtain the linear model for the voltage magnitudes $|\bv|$ in \eqref{eqn:lin_rho}, we leverage the following derivation rule: 
\[
\frac{\partial |f(x)|}{\partial x} = \frac{1}{|f(x)|} \Re\left\{ \conj{f(x)} \frac{\partial f(x)}{\partial x} \right\}.
\]
It then follows that matrices $\bK^Y$ and $\bK^\Delta$ are given by:
\begin{subequations} \label{eqn:K}
\begin{align}
\bK^Y := \frac{\partial |\bv|}{ \partial \bx^Y} = \diag (|\hat{\bv}|)^{-1}\Re \left \{ \diag (\conj{\hat{\bv}})  \bM^Y \right\}, \\
\bK^\Delta := \frac{\partial |\bv|}{ \partial \bx^\Delta} = \diag (|\hat{\bv}|)^{-1}\Re \left \{ \diag (\conj{\hat{\bv}})  \bM^\Delta \right\}, \\
\bb := |\hat{\bv}| - \bK^Y \hat{\bx}^Y - \bK^\Delta \hat{\bx}^\Delta.
\end{align}
\end{subequations}

\iffalse
Finally, we obtain \eqref{eqn:lin_s0} by using \eqref{eqn:lin_v} and \eqref{eqn:s0}, so that
\begin{subequations} \label{eqn:G}
\begin{align}
\bG^Y = \diag(\bv_0) \conj{\bY}_{0L} \conj{\bM}^Y, \, \bG^\Delta = \diag(\bv_0) \conj{\bY}_{0L} \conj{\bM}^\Delta, \\  
\bc = \diag(\bv_0) \left(\conj{\bY}_{00} \conj{\bv}_0 + \conj{\bY}_{0L} \conj{\ba}\right).
\end{align}
\end{subequations}
\fi

\subsection{Fixed-Point Linearization (FPL) Method}

Let $\hat{\bv}, \hat{\bs} := ((\hat{\bs}^Y)^\sfT, (\hat{\bs}^\Delta)^\sfT)^\sfT$ be a given solution to the fixed point equation \eqref{eqn:fp}. 
\rev{For a given power injection vector $\bs := ((\bs^{Y})^\sfT, (\bs^{\Delta})^\sfT)^\sfT$,}
consider the first iteration of the fixed-point method \eqref{eqn:fp_iter} initialized at $\hat{\bv}$:
\begin{equation} \label{eqn:lin_fp}
\tilde{\bv} = \bw + \bY_{LL}^{-1}\left(\diag(\conj{\hat{\bv}})^{-1}\conj{\bs^{Y}} + \bH^\sfT \diag\left(\bH \conj{\hat{\bv}} \right)^{-1} \conj{\bs^{\Delta}} \right)
\end{equation}
which gives an \emph{explicit} linear model \eqref{eqn:lin_v} provided by
\begin{align*}
&\bM^Y := \left(\bY_{LL}^{-1}\diag(\conj{\hat{\bv}})^{-1}, -\jay \bY_{LL}^{-1}\diag(\conj{\hat{\bv}})^{-1} \right) \\
&\bM^\Delta := \left(\bY_{LL}^{-1}\bH^\sfT \diag\left(\bH \conj{\hat{\bv}} \right)^{-1}, -\jay \bY_{LL}^{-1}\bH^\sfT \diag\left(\bH \conj{\hat{\bv}} \right)^{-1} \right)
\end{align*}
and $\ba = \bw$. The model \eqref{eqn:lin_rho} %and \eqref{eqn:lin_s0} 
can be then obtained \revvv{by substituting the above expressions for $\bM^Y$ and $\bM^\Delta$ in \eqref{eqn:K}}. %and \eqref{eqn:G}, respectively.
We next provide an upper bound for the linearization error of  the FPL method.
\begin{theorem} \label{thm:fixedAccuracy}
Suppose that $(\hat{\bv}, \hat{\bs})$ satisfy condition \eqref{eqn:cond1}. Let $\bs$ be the vector of power injections that satisfies \eqref{eqn:cond2}, and let $\bv \in \cD_\rho(\hat{\bv})$ with $\rho=\rho^\dagger(\hat{\bv},\hat{\bs},\bs)$ be the corresponding unique load-flow solution as guaranteed by Theorem \ref{thm:exist_lf}. Then the approximation error of \eqref{eqn:lin_fp} can be upper bounded by
\begin{equation}
\|\tilde{\bv} - \bv \|_\infty \leq q \rho^\dagger(\hat{\bv},\hat{\bs},\bs)\|\bw\|_\infty
\end{equation}
where
\[
q := \frac{\xi^Y(\mathbf{s})}{(\alpha(\hat{\bv})-\rho^\dagger(\hat{\bv},\hat{\bs},\bs))^2}+\frac{\xi^\Delta(\mathbf{s})}{(\beta(\hat{\bv})-\rho^\dagger(\hat{\bv},\hat{\bs},\bs))^2} < 1.
\]
\end{theorem}

\begin{figure}[t]
  \centering
  \includegraphics[width=0.9\columnwidth]{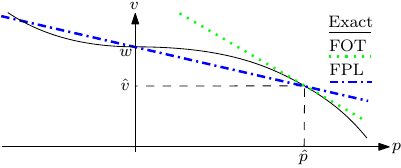}
  \vspace{-.3cm}
\caption{Qualitative interpretation of the FOT and FPL  methods.} \label{fig:illustr}
\vspace{-.5cm}
\end{figure}

The difference between the two linearization methods is conceptually illustrated in Figure \ref{fig:illustr}. The fixed-point linearization method can be viewed as an \emph{interpolation} method between two load-flow solutions: $(\bw, {\bf 0})$ and $(\hat{\bv}, \hat{\bs})$. On the other hand, the FOT yields the tangent plane of the load-flow manifold at the current linearization point.   

Some qualitative comparison between the FOT and FPL methods follows (a numerical comparison is provided shortly in Section \ref{sec:numerical}).
The FOT method provides the best local linear approximator, and hence it is expected to provide the best approximation accuracy around the linearization point. However, the main downside of the FOT method is its computational complexity. Indeed, solving $(12N)^2$ equations with $(12N)^2$ variables might not be feasible for large $N$ (i.e., large networks). On the other hand, the FPL method is computationally affordable as it requires only elementary vector-matrix multiplications (provided that $\bY_{LL}^{-1}$ is precomputed in advance). Moreover, if global behaviour is of interest, it can also provide a better approximation (cf.~Figure \ref{fig:illustr}). \rev{As a result, the FOT method may be preferable in a slowly time-varying setting whereby the variation of the power injections is relatively small. On the other hand, in the setting of modern distribution networks with high penetration of renewables, the FPL method may be preferable.}
\rev{
\begin{remark}
Using methods similar to the previous remarks, the results presented in this section can be straightforwardly adapted to the cases of general multiphase networks and the cases where no phase-to-phase connection exists.
\end{remark}
}
%Interestingly, we will show in Section ?? that it may provide a better approximation than the FOT method in the reverse-power-flow scenarios which is typical in modern distribution networks.     

%One of the advantages of using \eqref{eqn:lin_fp} is that it represents a single iteration of the the fixed-point equation \eqref{eqn:fp}. So if this model is utilized in an online control algorithm, it can have nice convergence properties (it is a conjecture, for future research...)
%Another advantage is that the coefficient matrix can be updated locally, thus allowing for completely distributed implementation of controllers. \andy{Give more details on that.}

\rev{
\section{Potential Applications}
In this section, we briefly discuss the potential applications of our results.  As mentioned in the introduction, they can be used to facilitate the development of OPF solvers and  real-time control procedures for general multiphase distribution networks. In particular:
\begin{itemize}
\item  Linear models of Section \ref{sec:linear} can be leveraged to convexify the OPF problem, and thus facilitate the development of OPF-based real-time control techniques.
Particularly, the  methodology proposed in this paper can be utilized to broaden the applicability of~\cite{swaroop2015linear,commelec1,opfPursuit,voltVARWattArxiv} to the case of unbalanced multiphase systems with delta and wye connections. 
\item Explicit conditions of Theorem \ref{thm:exist_lf} can be directly embedded in the optimization problems as convex constraints, thus ensuring existence and non-singularity of the exact high-voltage load-flow solution. 
\end{itemize}
}

\section{Numerical Evaluation} \label{sec:numerical}
\rev{In this section, we evaluate numerically the proposed methodology using IEEE test feeders \cite{testfeeder}. Particularly, in the IEEE 37-Bus and 123-Bus networks, we compare our method with the method in \cite{ZIP}, which is the classic $Z$-bus method applied to the multiphase setting with disjoint sets of wye- and delta-connected sources. 
We also use the IEEE 8500-Node test feeder to demonstrate the applicability of the proposed algorithms to a large-scale distribution network.
%From the first to the last example, we will gradually show networks of different sizes, sources/loads connections, number of phases, etc. Also, we will give further explanation of some important concepts proposed in this paper.
}

\begin{figure}[t]
  \centering
  \includegraphics[scale=0.33]{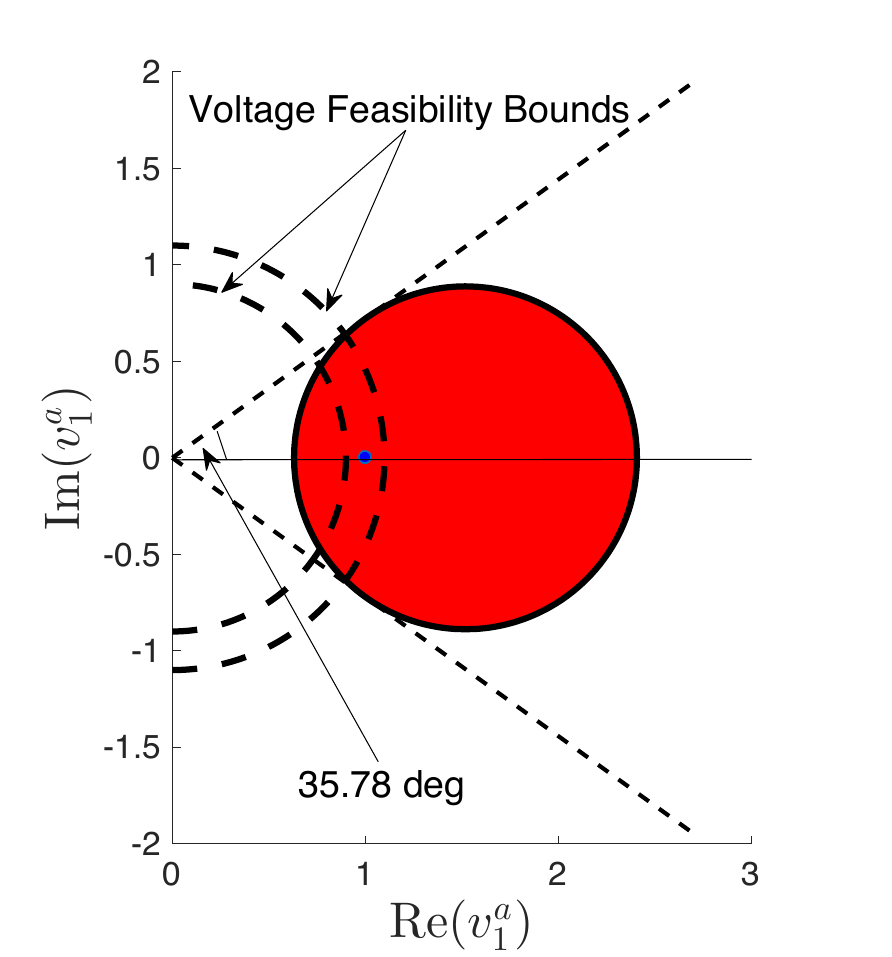}
  \includegraphics[scale=0.33]{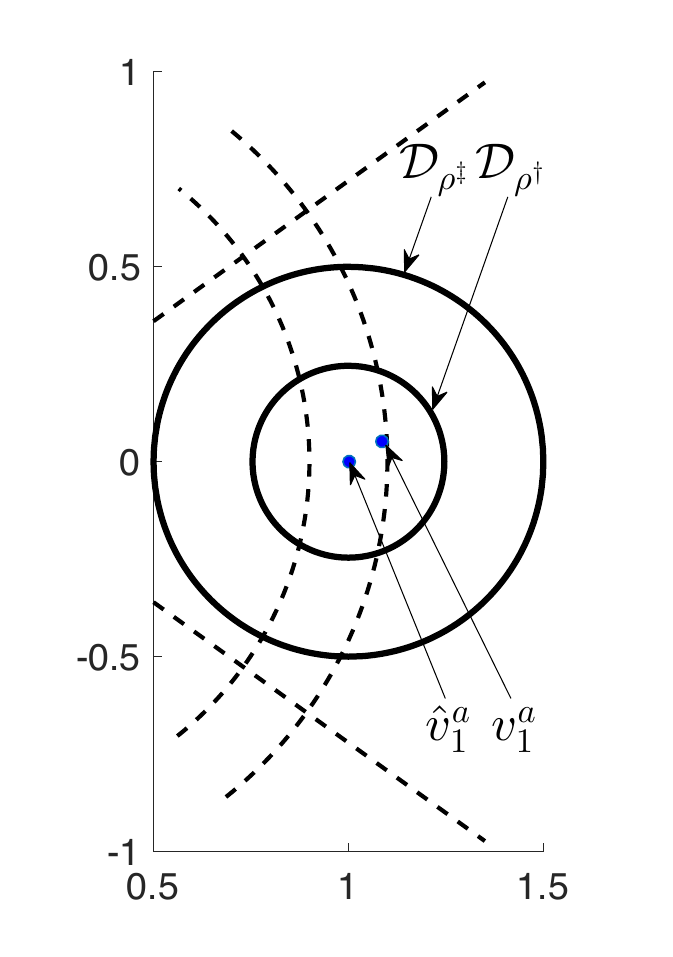}
\caption{\rev{Illustration in voltage space, where the unit of axes is p.u. (Left) The \revv{filled red circle represents the set} where condition \eqref{eqn:cond1} holds, \revv{which contains the zero-load point $v_1^a=1$ (shown as a blue dot)}; (Right) The set $\mathcal{D}_\rho$ \revv{(projected onto the space of $v_1^a$)} and the guaranteed solution.}} \label{fig:SimpleExampleCondD}
\end{figure}

\rev{
\subsection{An Illustrative Example} \label{sec:simpleillu}
We start by demonstrating the proposed methodology  and its physical significance using an artificially-designed network. Here, the purpose is to facilitate the understanding and meanwhile provide some intuition.

We consider a balanced network with a single three-phase $PQ$ bus (with index 1) connected to the slack bus via a transmission line. The line admittance matrix is given as follows, in p.u.:
\begin{equation}
\begin{bmatrix}
7-12\jmath & -1+2\jmath & -1+2\jmath \\
-1+2\jmath & 7-12\jmath & -1+2\jmath \\
-1+2\jmath & -1+2\jmath & 7-12\jmath
\end{bmatrix}.
\end{equation}
Moreover, we assume that the shunt elements are negligible and the vector of slack-bus voltages is $\bv_0=(1,e^{-\jmath\frac{2\pi}{3}},e^{\jmath\frac{2\pi}{3}})^T$ p.u. Therefore,
\begin{equation}
\bw=\bv_0,~\bY_{LL}=
\begin{bmatrix}
7-12\jmath & -1+2\jmath & -1+2\jmath \\
-1+2\jmath & 7-12\jmath & -1+2\jmath \\
-1+2\jmath & -1+2\jmath & 7-12\jmath
\end{bmatrix}.
\end{equation}
Now, exclude the delta connections and let the power injection vector $\bs^Y$ be balanced in all phases. As a direct consequence, $\bv_1=(v_1^a)\bv_0$, which means that the vector of voltages at bus 1 is determined by a scalar $v_1^a$.

In the left-hand side of Figure \ref{fig:SimpleExampleCondD}, we plot the region (a filled circle) in the voltage space where condition \eqref{eqn:cond1} holds. It can be seen that this region covers almost all the $v_1^a$ with a feasible magnitude and an angle between $\pm 35.78^{\circ}$, which is of practical significance. Also, note that the region contains $v_1^a$ with a magnitude much higher than 1 p.u., which corresponds to the case of strong reverse power flow. In the right-hand side of Figure \ref{fig:SimpleExampleCondD}, we take $\hat{v}_1^a=1$ p.u. \revv{(i.e., $\hat{\bv}=\bw$, $\hat{\bs}=\mathbf{0}$)}, $\bs^Y=(1.5+0.9\jmath,1.5+0.9\jmath,1.5+0.9\jmath)^T$ p.u., and plot the domain $\mathcal{D}_\rho(\hat{\bv})$ projected on $v_1^a$ for the typical radii \revv{$\rho^\ddagger(\hat{\bv},\hat{\bs})$, $\rho^\dagger(\hat{\bv},\hat{\bs},\bs)$} in \eqref{eqn:rhos}. We also show the solution $v_1^a$ in $\mathcal{D}_\rho(\hat{\bv})$ with $\rho=\rho^\ddagger(\hat{\bv},\hat{\bs})$, where $\hat{\bs}$ is the power injection corresponding to $\hat{\bv}$. It can be seen that, when taking the power injections vector $\bs$ into account, the guaranteed solution is localized more accurately using $\mathcal{D}_\rho(\hat{\bv})$ with $\rho=\rho^\dagger(\hat{\bv},\hat{\bs},\bs)$. In Table \ref{table:Tab0}, we present the update of ${v_1^a}^{(k)}$ during the iteration. By observing the third column, it is clear that the iterative update gradually converges. In the fourth column, we give the convergence rate, which is bounded by the contraction modulus  $\frac{\xi(\bs)}{\left(\gamma(\hat{\bv})-\rho^\dagger(\hat{\bv},\hat{\bs},\bs)\right)^2}=0.3264$ (see Appendix for reference).

We note that empirical evidences show that the true convergence rate is usually less than a third of the contraction modulus. As a consequence, when our conditions hold, the iterative method generally reaches a precision of $10^{-6}$ in less than ten iterations.
}

\begin{table}[t]
\caption{\rev{Update of $v_1^a$ in the iteration. In this table, values are presented with four decimal digits.}}
\label{table:Tab0}
\begin{center}
\rev{
\begin{tabular}{cccc}
\noalign{\global\arrayrulewidth0.05cm}
\hline
\noalign{\global\arrayrulewidth0.4pt}
$k$ & ${v_1^a}^{(k)}$ & $|{v_1^a}^{(k)}-{v_1^a}^{(k-1)}|$ & $\frac{|{v_1^a}^{(k+1)}-{v_1^a}^{(k)}|}{|{v_1^a}^{(k)}-{v_1^a}^{(k-1)}|}$ \\
\hline
$0$ & $1.0000+0.0000\jmath$ &  &  \\
$1$ & $1.0946+0.0531\jmath$ & 0.1085 & 0.0990 \\
$2$ & $1.0839+0.0526\jmath$ & 0.0107 & 0.0912 \\
$3$ & $1.0847+0.0531\jmath$ & 0.0010 & 0.0921 \\
$4$ & $1.0846+0.0531\jmath$ & 0.0001 &  \\
\noalign{\global\arrayrulewidth0.05cm}
\hline
\noalign{\global\arrayrulewidth0.4pt}
\end{tabular}
}
\end{center}
\end{table}

\subsection{IEEE 37-Bus Feeder}
\rev{In this example, we evaluate the performance of our method on a network with purely delta connections.} Similar to prior works~\cite{Low-Convex,swaroop2015linear,commelec1,opfPursuit,bolognani2015linear}, we translate all constant-current and constant-impedance sources in the IEEE data set into constant-power sources. \rev{In addition, we fix the voltage regulators in this and all subsequent examples at their default values.} %and remove %the transformer between buses 709 and 
%bus 775 (no loads are present). 

\begin{figure}[t!]
     \centering
     \begin{subfigure}[t]{0.49\textwidth}
         \centering
 		 \includegraphics[width=0.9\columnwidth]{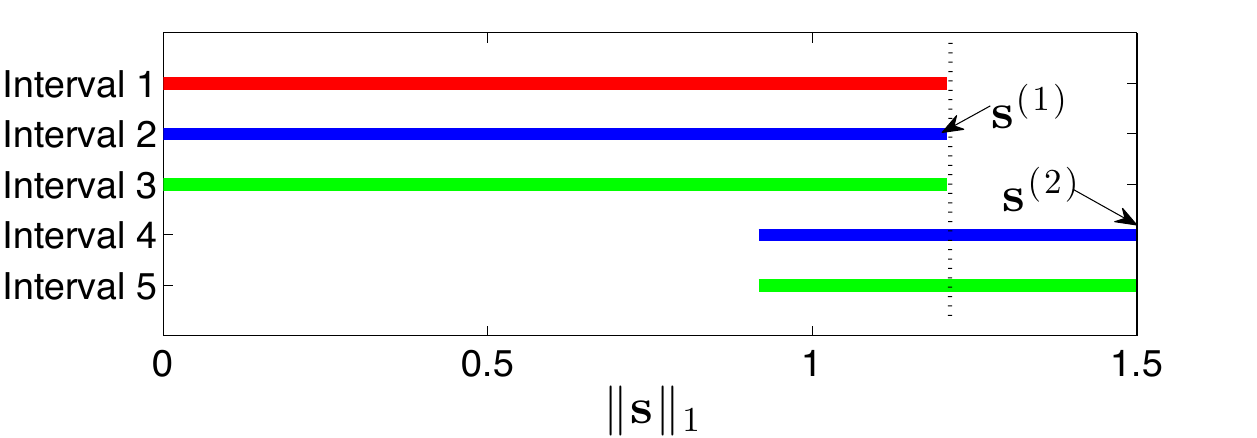}
		 \caption{\footnotesize Power intervals (in p.u.) that satisfy conditions on the existence and uniqueness of the load-flow solutions (Intervals 2, 3, 4, and 5) and comparison to the interval obtained by 	\cite{ZIP} (Interval 1).} \label{fig:37nodeinterval1}
     \end{subfigure} \\
     
     \begin{subfigure}[b]{0.49\textwidth}
         \centering
  		 \includegraphics[width=0.75\columnwidth]{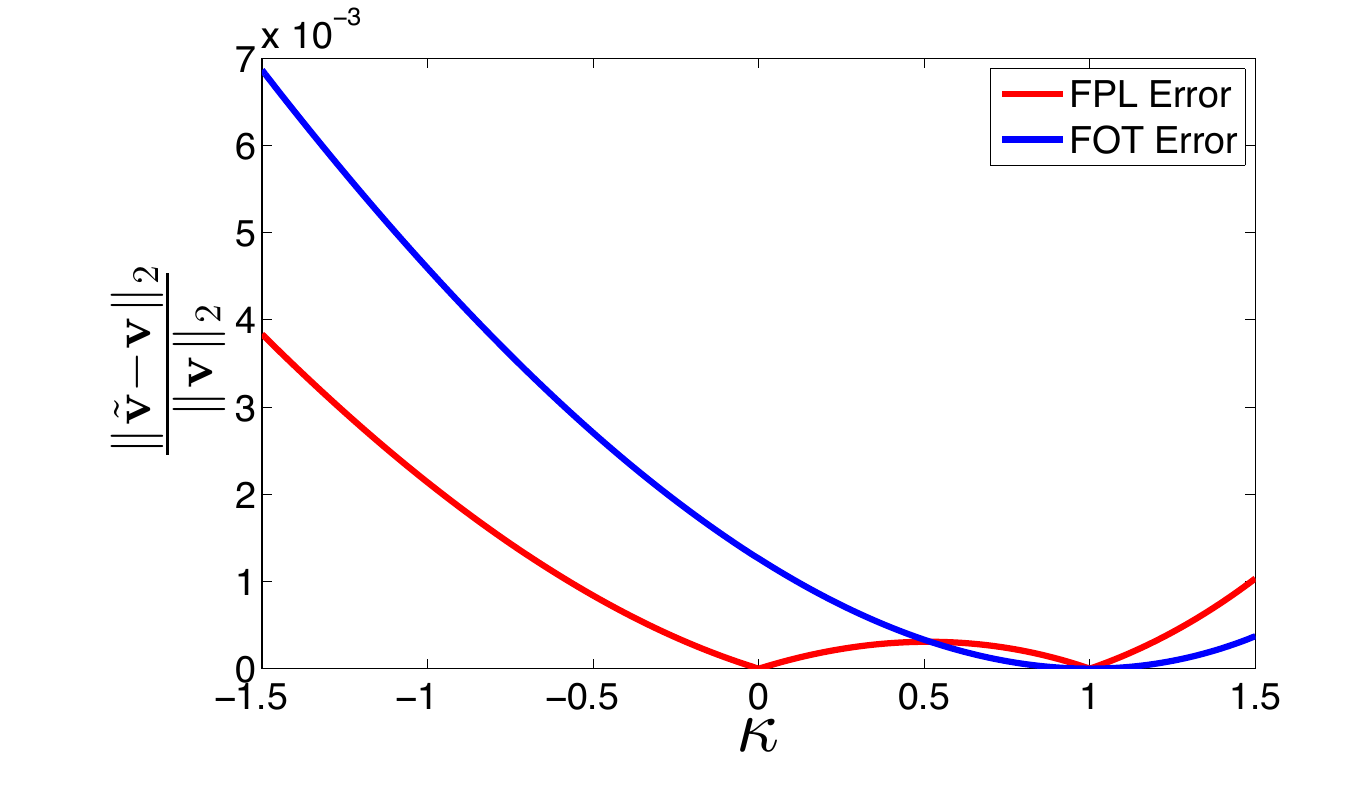}
		 \caption{\footnotesize Relative errors in complex load-flow solutions for FOT and FPL methods.} \label{fig:37nodelinear}
     \end{subfigure}
     \caption{Numerical evaluation for the 37-bus feeder.}
\end{figure}

In the original IEEE data set, all sources/loads are delta-connected. 
Denote this reference power injection vector by $\mathbf{s}^\mathrm{ref}$, and let the target power injection be $\mathbf{s}=\kappa\mathbf{s}^\mathrm{ref}$ with $\kappa$ as a real number.
As there are no mixed wye and delta sources/loads, the conditions on the existence and uniqueness of the load-flow solutions in \cite{ZIP} are also applicable. For comparison, we take the diagonal matrix $\mathbf{\Lambda}$ in \cite{ZIP} to be $\mathbf{W}$, as suggested there. In  Figure \ref{fig:37nodeinterval1}, we let $\kappa$ be nonnegative and plot five power intervals. in p.u. %In terms of the first three, we have:
Interval 1 contains the power injection $\mathbf{s}$ that satisfies the four conditions in \cite{ZIP}; Interval 2 (resp.~3) shows the injections $\mathbf{s}$ that satisfy the conditions in Theorem \ref{thm:exist_general} (resp.~Theorem \ref{thm:exist_lf}) with $(\hat{\bv},\hat{\bs})=(\bw,\mathbf{0})$. For the rightmost power $\mathbf{s}^{(1)}=3.45\mathbf{s}^\mathrm{ref}$, we  compute the load-flow solution $\bv^{(1)}$ using iteration \eqref{eqn:fp_iter} \rev{(initialized at $\bw$)}. By choosing this solution $\bv^{(1)}$ and $\mathbf{s}^{(1)}$ as the new $(\hat{\mathbf{v}},\hat{\mathbf{s}})$, we obtain Interval 4 (resp. 5) via Theorem \ref{thm:exist_general} (resp.~Theorem \ref{thm:exist_lf}). \rev{Note that for this choice of $(\hat{\bv},\hat{\bs})$, only some of the power injections in Interval 2 (resp. 3) satisfies the proposed conditions. This is because the conditions guarantee the solution properties only for the power injections in a domain around $\hat{\bs}$. It can be further shown that, for any power injection vector $\bs$ in the intersection of Interval 2 (resp. 3) and Interval 4 (resp. 5), the guaranteed 
load-flow solution $\bv$ is consistent. This is because $\bv$ can be computed by iteration \eqref{eqn:fp_iter} initialized at $\bv^{(1)}$.}

Numerically, Intervals 1,2, and 3 are the same. However, the complexity of computing  Interval 3 is much smaller because of  the low computational complexity of verifying conditions \eqref{eqn:cond1} and \eqref{eqn:cond2}.
More importantly, Intervals 4 and 5  
%In addition, point like $\mathbf{s}^{(2)}=4.28\mathbf{s}^\mathrm{ref}$ 
contain points that are not guaranteed to have the unique solution using the method in \cite{ZIP} -- compare to Interval 1. Thus, the proposed method allows for certifying the existence and uniqueness of the load-flow solution for a wider range of power injections.  

Next, we evaluate the performance of the two linearization methods proposed in Section \ref{sec:linear}. 
Figure \ref{fig:37nodelinear} shows the results of the relative errors for both linear models using $\kappa\in[-1.5,1.5]$. 
%Particularly, for the injections vector $\mathbf{s}^{(2)}=1.60\mathbf{s}^\mathrm{ref}$ shown in Figure \ref{fig:37nodeinterval2}, we compute the relative error between the solution provided by the FOT and FPL methods with respect to the exact solution.
As shown, both linear models behave well with relative errors below $1\%$. Moreover, the FOT method has a smaller error around the linearization point whereas the FPL method provides a better global approximation.  This corroborates the intuitive illustration in Figure \ref{fig:illustr}. For linear approximations of voltage magnitudes, %and power injections through the slack bus, 
the errors are at a similar level; hence, for brevity, we do not show them explicitly. %In the next section, we will present an evaluation of the linear models' performance  on a more global scale.
%\andy{...}

\iffalse
Based on $\hat{\mathbf{v}},\hat{\mathbf{s}}$, we can compute the solution to power injection $\mathbf{s}^{(2)}=1.60\mathbf{s}^\mathrm{ref}$ either exactly via iteration \eqref{eqn:fp_iter} or approximately via FOT and FPL methods given in Section \ref{sec:linear}. As it is of interest to have a glimpse of how linear models behave, we compute the relative errors for complex load-flow solution and depict them in Figure \ref{fig:37nodelinear}. It can be seen that both linear models behave well with bus-wise relative errors lower than $0.3\%$. Moreover, FOT method has a smaller error, which in some sense mirrors the intuitive illustration in Figure \ref{fig:illustr}.
\fi

\vspace*{-0.3cm}
\subsection{IEEE 123-Bus Feeder} \label{sec:feed123}
In this section, we consider a larger multiphase network with unbalanced one-, two-, and three-phase sources/loads. This network represents the normal size of many distribution networks in the world. As mentioned in Remark \ref{rem:missingPhases}, we first delete in matrix $\bH$ the rows that correspond to the lacking phase-to-phase connections and the columns that correspond to the lacking phases.

%Unlike the previous example, this network is larger in size and lacks some phases in certain buses. Here, we will first demonstrate that the proposed theory works in this larger multiphase network, and then complete the missing part of linear model analysis. Before proceeding, as mentioned in Remark \ref{rem:missingPhases}, we should delete in $\bH$ the rows that correspond to non-existent between-phase connections and the columns that correspond to non-existence phases.
\iffalse
\begin{figure}[t]
  \centering
  \includegraphics[width=0.9\columnwidth]{123NodeIntervalNew}
\caption{Conditions evaluation for 123-bus feeder using original injections data. Power intervals (in p.u.) that satisfy conditions on existence and uniqueness of load-flow solutions (Intervals 2, 3, 4, and 5) and comparison to the interval obtained by \cite{ZIP} (Interval 1).} \label{fig:123nodeinterval1}
\vspace{-.5cm}
\end{figure}

\begin{figure}[t]
  \centering
  \includegraphics[width=0.9\columnwidth]{123NodeIntervalNonTrivialNew}
\caption{Conditions evaluation for 123-bus feeder with mixed delta-wye sources. Power intervals (in p.u.) that satisfy conditions on existence and uniqueness of load-flow solutions.} \label{fig:123nodeinterval2}
\vspace{-.5cm}
\end{figure}
\fi

\begin{figure}[t!]
     \centering
     \begin{subfigure}[t]{0.49\textwidth}
         \centering
  		\includegraphics[width=0.9\columnwidth]{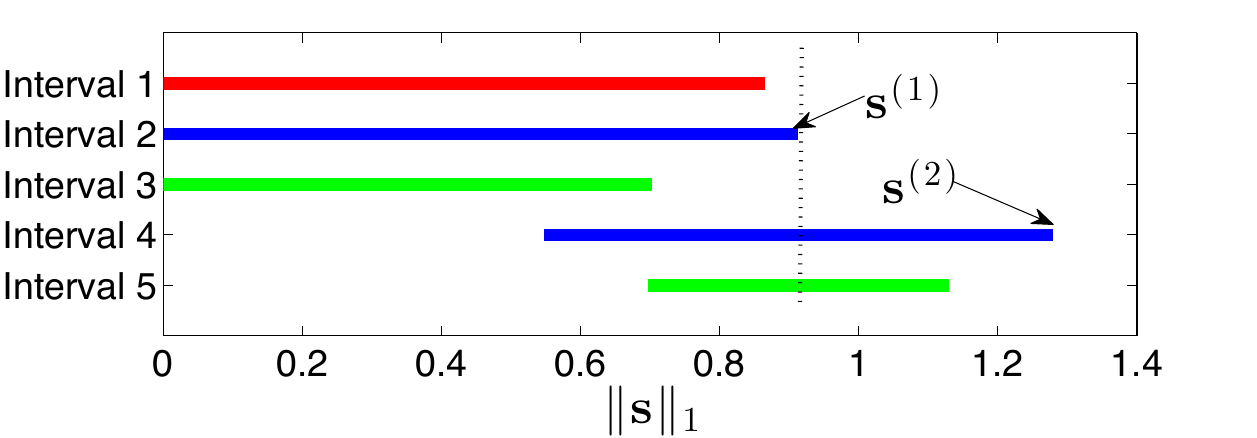}
		\caption{\footnotesize Original injections data. Power intervals (in p.u.) that satisfy conditions on the existence and uniqueness of the load-flow solutions (Intervals 2, 3, 4, and 5) and comparison to the interval obtained by \cite{ZIP} (Interval 1).} \label{fig:123nodeinterval1}
     \end{subfigure} \\
     \begin{subfigure}[b]{0.49\textwidth}
         \centering
  		 \includegraphics[width=0.9\columnwidth]{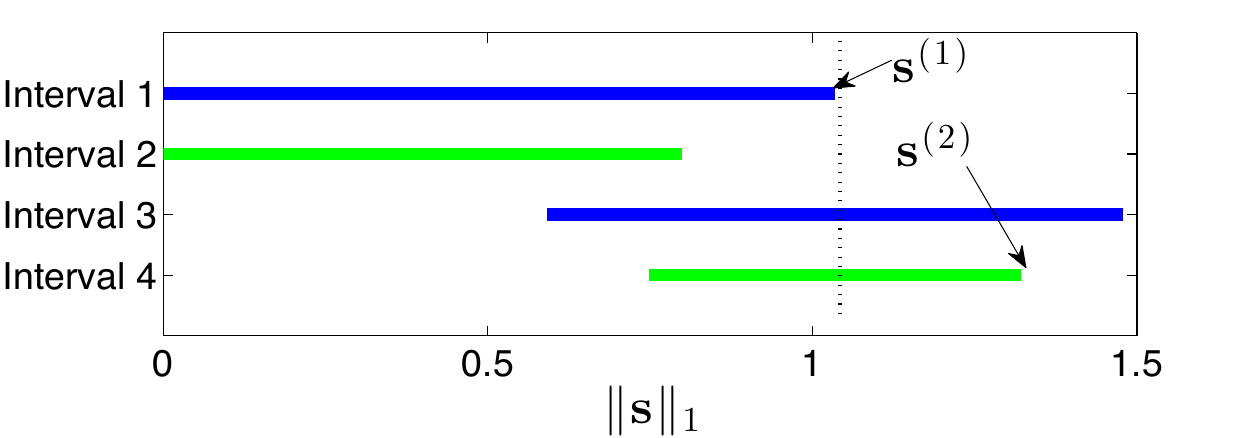}
		 \caption{\footnotesize Mixed delta and wye sources. Power intervals (in p.u.) that satisfy conditions on the existence and uniqueness of the load-flow solutions.} \label{fig:123nodeinterval2}
     \end{subfigure}
     \begin{subfigure}[b]{0.49\textwidth}
         \centering
  		\includegraphics[width=0.75\columnwidth]{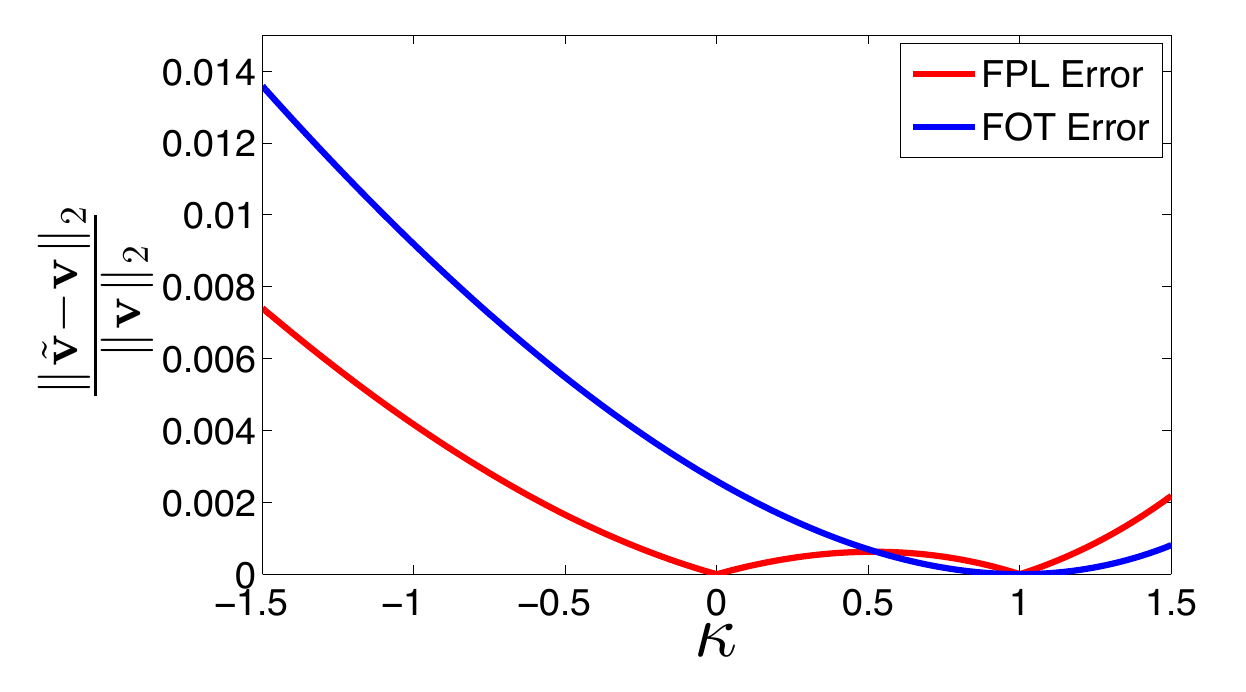}
		\caption{\footnotesize Relative errors in complex load-flow solutions for FOT and FPL methods.} \label{fig:123nodelinear}
     \end{subfigure}
     \caption{Conditions evaluation for the 123-bus feeder.}
\end{figure}

\rev{Similar to the previous case, let $\bs=\kappa\bs^\mathrm{ref}$ with $\bs^\mathrm{ref}$ being the reference power injections in this network. Consider then  repeating the analysis of the previous subsection. The results are shown in Figures \ref{fig:123nodeinterval1}, with the same interpretation of the intervals as in Figures \ref{fig:37nodeinterval1}. To perform the experiment with mixed delta-wye connections, additional power sources/loads were added to the network, as shown in  Table \ref{table:Tab2}.  In this case of mixed connections, we obtain the intervals of \ref{fig:123nodeinterval2} in a way similar to the previous analysis. The results match with those obtained in OpenDSS \cite{openDSS}, which is the only freely-available solver that works with mixed connections.}

\iffalse
In the IEEE data set, each of the given sources are either wye-conneted or delta-connected. In the foremost place, we apply Theorem \ref{thm:exist_general} and confirm that this benchmark power injection satisfies out proposed conditions with $(\hat{\bv},\hat{\bs})=(\bw,\mathbf{0})$. Then, we use the iterative method (initialized at $\bw$) to obtain the load-flow solution for this benchmark power injection, and compare it with the one obtained from OpenDSS. The two results match exactly, which means that the proposed theory can handle multiphase network by simply adapting matrix $\bH$.

Then, we create some mixed delta-wye sources according to Table \ref{table:Tab2}. Like before, denote the IEEE benchmark power injection together with these additional power sources by $\mathbf{s}^\mathrm{ref}$, and let the target power injection be $\mathbf{s}=\kappa\mathbf{s}^\mathrm{ref}$ with $\kappa$ being a real number. Note that we allow reverse power flow (i.e., cases with negative $\kappa$) in this example, which is permitted by our usage of $\ell_\infty$ norms and absolute values.
\fi

\begin{table}[t]
\caption{Additional Power Sources in 123-Bus Feeder}
\label{table:Tab2}
\begin{center}
\begin{tabular}{ccccc}
\noalign{\global\arrayrulewidth0.05cm}
\hline
\noalign{\global\arrayrulewidth0.4pt}
Bus & Type & Phase-Phase ab & Phase-Phase bc & Phase-Phase ca \\
& & / Phase a (p.u.) & / Phase b (p.u.) & / Phase c (p.u.) \\
\hline
$1$ & delta & -0.03-0.01$\jmath$ & -0.03-0.01$\jmath$ & -0.03-0.01$\jmath$ \\
$35$ & wye & -0.02 & -0.02 & -0.02 \\
$76$ & wye & 0.04+0.01$\jmath$ & 0.04+0.01$\jmath$ & 0.04+0.01$\jmath$ \\
$99$ & delta & -0.02-0.01$\jmath$ & -0.02-0.01$\jmath$ & -0.02-0.01$\jmath$ \\
\noalign{\global\arrayrulewidth0.05cm}
\hline
\noalign{\global\arrayrulewidth0.4pt}
\end{tabular}
\end{center}
\vspace{-.5cm}
\end{table}

Finally, in Figure \ref{fig:123nodelinear}, we show the results of the relative errors for both linear models using $\kappa\in[-1.5,1.5]$. \rev{Here, different from the counterpart in the last section, we have incorporated in $\bs^\mathrm{ref}$ the additional sources in Table \ref{table:Tab2}.} Clearly, the errors vary in a way that is similar to the illustration in Figure \ref{fig:illustr}. In other words, the FPL method provides not only a high computational efficiency but also a better global performance for large distribution networks.

%\rev{In summary, through this example, we demonstrate the validity of proposed theories in a typical multiphase network with mixed wye and delta sources/loads.}

\iffalse
By application of Theorem \ref{thm:exist_general} with $(\bw,\mathbf{0})$, it can be checked that all power injections with $\kappa\in[-1.05,1.05]$ satisfies the proposed conditions and hence have guaranteed unique solutions around $\bw$. Moreover, it can be verified that both linear models can be established at the solved load-flow solution for $\kappa=1$. In this way, for each $\kappa\in[-1.05,1.05]$, we have an exact solution obtained via iterations and two approximate solutions obtained via different linear models. In Figure \ref{fig:123nodelinear}, we plot the vector-wise relative errors for both linear models. Clearly, the errors vary in a way similar to the illustration in Figure \ref{fig:illustr}. In other words, FPL method provides not only a high computational efficiency but also better global performance; therefore, it is of practical interest.
\fi

\iffalse
\begin{figure}[t]
  \centering
  \includegraphics[width=0.9\columnwidth]{123LinearErrorExpand}
\caption{Continuation analysis for linear models in 123-bus feeder: Relative errors in complex load-flow solutions for FOT and FPL methods.} \label{fig:123nodelinear}
\vspace{-.5cm}
\end{figure}
\fi

\rev{
\subsection{IEEE 8500-Node Feeder}
In this subsection, we illustrate the performance of the proposed methodology using the IEEE 8500-Node feeder \cite{8500testFeeder}.  This network represents a large-scale distribution network with detailed modeling of the secondary side of distribution transformers.
%Same as before, we disable the regulators. 

In this network, the line-to-line medium-voltage rating is 12.47 kV, and the network contains split-phase secondary loading with  line-to-line low-voltage rating of 208 V. In Figure  \ref{fig:range8500}, we evaluate the working range of the proposed methodology. \revv{In particular, let $\bs=\kappa\bs^\mathrm{ref}$ and $\bv$ be the guaranteed load-flow solution that corresponds to $\bs$. Moreover, define the feasibility constraints as $|(\bv)_j|\geq 0.9|(\bw)_j|,\forall j$, where $\bw$ is the zero-load voltage profile given in \eqref{eqn:fp}. In this way, $\alpha(\bv)$ (defined in \eqref{eqn:gamma}) becomes both a function of $\|\bs\|_1$ and an indicator of the feasibility.

Now,} given the knowledge of the zero-load voltage $\bw$, the maximum (in terms of $\ell_1$-norm) power vector that satisfies conditions \eqref{eqn:cond1} and \eqref{eqn:cond2} is $\bs^{(1)}$. Since the conditions are satisfied, we solve for its load-flow solution $\bv^{(1)}$. From the figure, it can be seen that there is already some voltage close to the feasibility boundary. Next, we take the values of $\hat{\bv}$ (resp. $\hat{\bs}$) to be $\bv^{(1)}$ (resp. $\bs^{(1)}$). Applying again the proposed conditions, we obtain that the maximum power vector is $\bs^{(2)}$, and the corresponding load-flow solution $\bv^{(2)}$ is obtained. As shown in the figure, some of the voltages in vector $\bv^{(2)}$ are already out of the feasibility region. By taking $\hat{\bv}$ (resp. $\hat{\bs}$) to be $\bv^{(2)}$ (resp. $\bs^{(2)}$), we continue the above procedure. Clearly, for this network, some of the voltages drop quickly due to its configuration and the disabled voltage regulators. 
Because our conditions rely on the voltages, their application becomes more challenging; however, we demonstrate that the conditions can be applied even in the cases where the voltages are significantly below the voltage feasibility boundary.

Finally, in Figure \ref{fig:8500_linear}, we evaluate the performance of the FPL method for this test feeder. Specifically, we plot the relative error of the phasor approximation using \eqref{eqn:lin_fp} and the corresponding magnitudes approximation using \eqref{eqn:K}, for $\kappa\in[-1,2]$. It can be seen that the relative errors are below $1.4\%$, confirming good scalability of the proposed linear approximation methodology for large-scale distribution networks.

\begin{figure}[t]
  \centering
  \includegraphics[width=0.9\columnwidth]{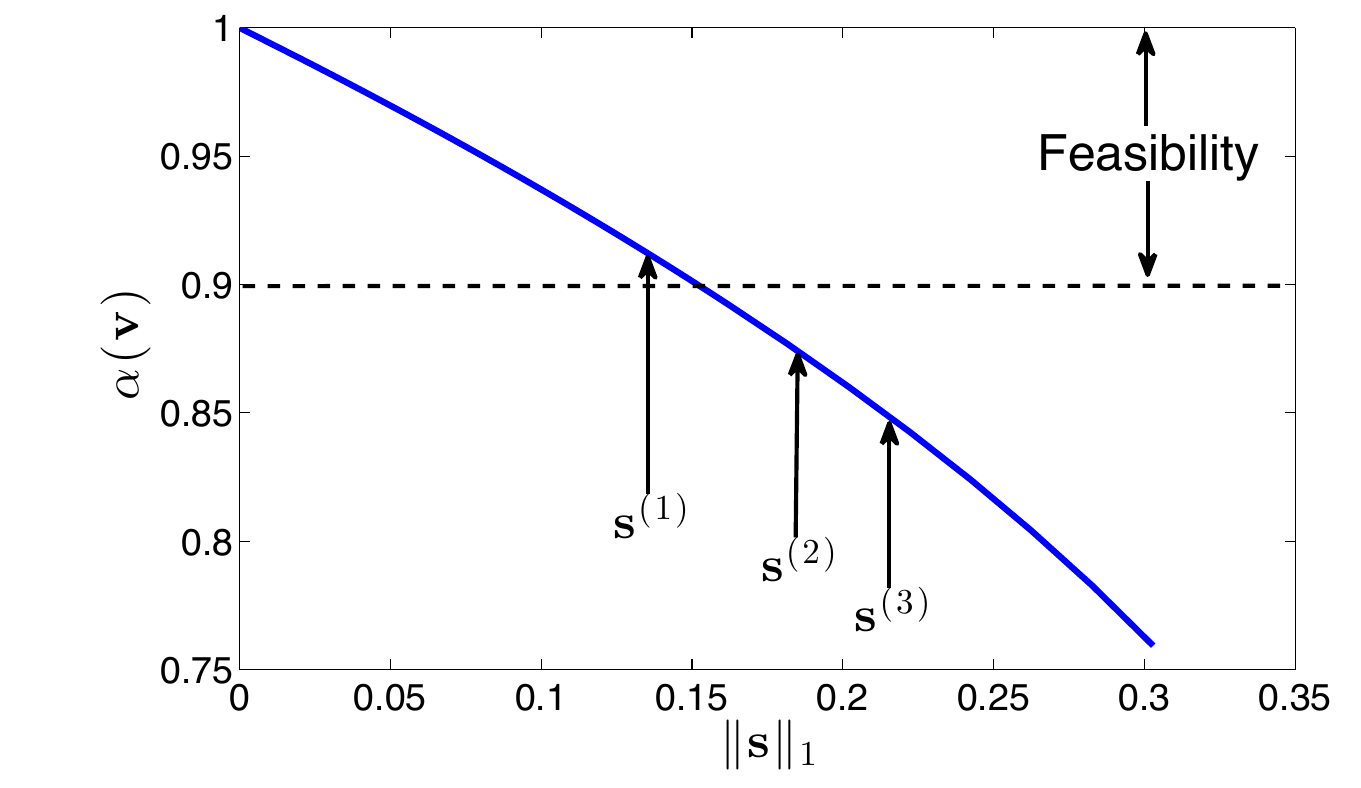}
\caption{\footnotesize 8500-node test feeder: illustration of the applicability of conditions.} \label{fig:range8500}
\end{figure}

\begin{figure}[t]
  \centering
  \includegraphics[width=0.75\columnwidth]{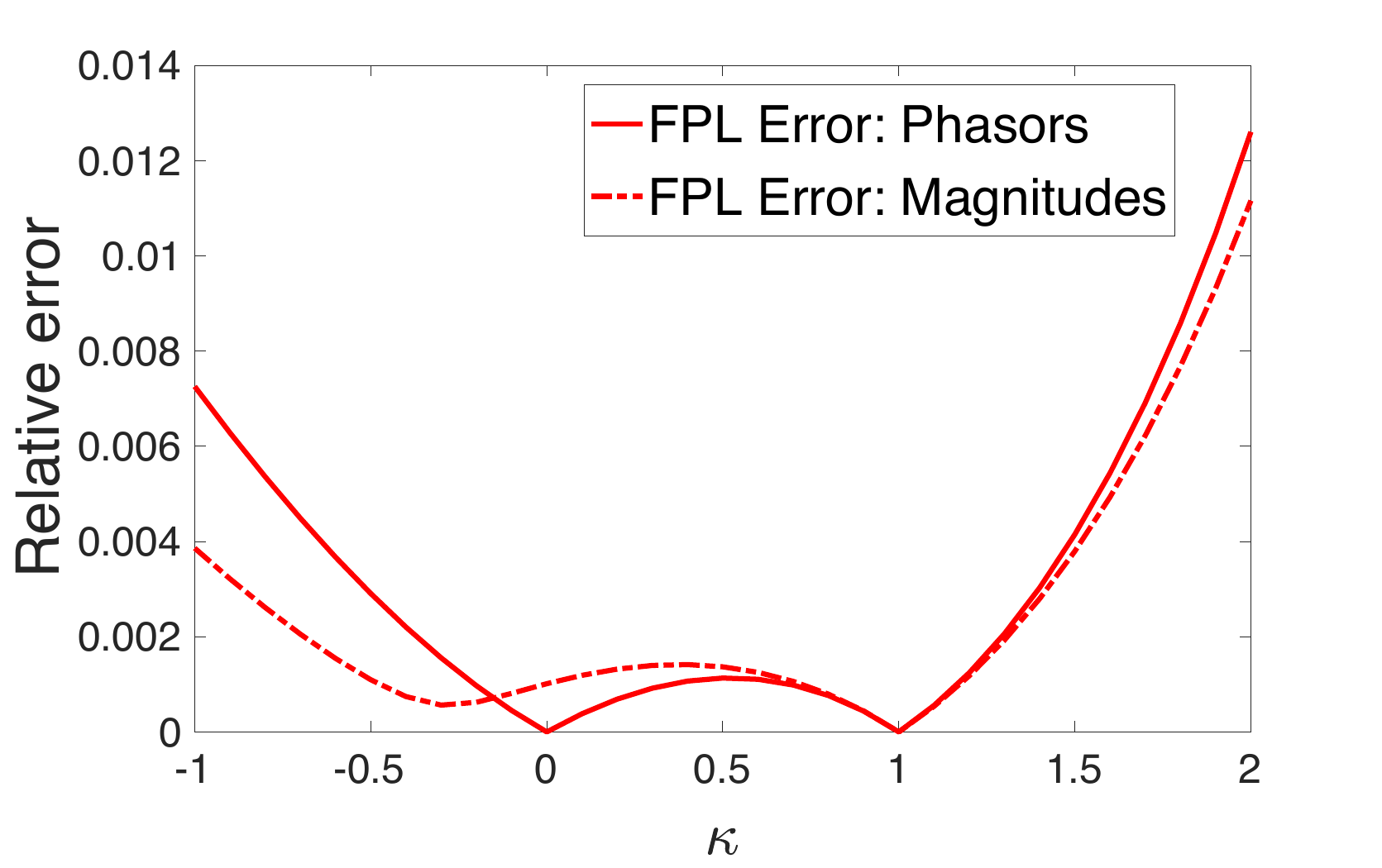}
\caption{\footnotesize 8500-node test feeder: relative errors in complex load-flow solutions and the corresponding magnitudes approximation for the FPL method.}  \label{fig:8500_linear}
\end{figure}

\subsection{Complexity Evaluation}
We next analyze the computational complexity of the proposed algorithms. In particular: (i) the verification of conditions \eqref{eqn:cond1} and \eqref{eqn:cond2} mainly depends on the computation of $\xi(\cdot)$ defined in \eqref{eqn:xi}, which has a worst-case complexity of $\mathcal{O}(N^2)$; (ii) the FPL linear model is essentially a single iteration of \eqref{eqn:fp_iter}, which takes $\mathcal{O}(N)$ to complete with LU decomposition in radial networks. To confirm the analysis, we measure the CPU time using MATLAB (on Macbook Pro @3GHz) and gather the results in Table \ref{table:Tab3}. From the second column of Table \ref{table:Tab3}, note that the conditions \eqref{eqn:cond1} and \eqref{eqn:cond2} can be verified efficiently for 37-Bus and 123-Bus networks, but cannot be verified in real-time for the 8500-Bus network. This adds some restrictiveness in the online applications to very large networks. However, when we pay attention to the third column, the complexity of the FPL method (i.e., single iteration of \eqref{eqn:fp_iter}) scales well with respect to the network size. Recall that, in almost all the experiments, the required number of iterations for accuracy $10^{-6}$ is less than 10. Therefore, the proposed methodology can be very useful in the real-time control and OPF in large networks.

\begin{table}[t]
\rev{
\caption{Complexity Evaluation}
\label{table:Tab3}
\begin{center}
\begin{tabular}{ccc}
\noalign{\global\arrayrulewidth0.05cm}
\hline
\noalign{\global\arrayrulewidth0.4pt}
Network & Condition \eqref{eqn:cond1} and \eqref{eqn:cond2} & Single iteration of \eqref{eqn:fp_iter} \\
\hline
$37$-Bus & $2.3$ ms & $0.17$ ms \\
$123$-Bus & $12$ ms & $0.49$ ms \\
$8500$-Bus & $76$ s & $51$ ms \\
\noalign{\global\arrayrulewidth0.05cm}
\hline
\noalign{\global\arrayrulewidth0.4pt}
\end{tabular}
\end{center}
}
\end{table}
}

% \begin{itemize}
% \item Show regions of convergence
% \item Compare the linearization methods to the exact solution
% \end{itemize}

% \subsection{Comparison with OpenDSS solution} \label{sec:compOpenDSS}

% \begin{figure}[t!]
%     \centering
%     \begin{subfigure}[t]{0.5\textwidth}
%         \centering
%         \includegraphics[width=1.0\columnwidth]{F_test_37_original}
%         \caption{Lorem ipsum}
%     \end{subfigure} \\
%     \begin{subfigure}[b]{0.5\textwidth}
%         \centering
%         \includegraphics[width=1.0\columnwidth]{F_test_37_modified}
%         \caption{Lorem ipsum}
%     \end{subfigure}
%     \begin{subfigure}[b]{0.5\textwidth}
%         \centering
%         \includegraphics[width=1.0\columnwidth]{F_conv}
%         \caption{Lorem ipsum}
%     \end{subfigure}
%     \caption{Caption place holder}
% \end{figure}

%\vspace*{-0.8cm}
\section{Conclusion} \label{sec:conclusion}

The paper \rev{extended the classical $Z$-bus} load-flow algorithm to general multiphase distribution systems. \rev{We derived explicit conditions for the existence of the load-flow solution, and analytically specified a domain in which the solution is unique. These conditions also guarantee the convergence of the load-flow algorithm to this solution. Then, we gave a sufficient condition for the non-singularity of the load-flow Jacobian, and proved that our theoretically guaranteed solution automatically ensures the non-singularity of the load-flow Jacobian. Finally, linear load-flow models were proposed and their approximation accuracy was analyzed.} Theoretical results were corroborated through numerical experiments on the IEEE test feeders. 

\rev{As we have discussed in the paper, the proposed theory and methodology can be leveraged in real-time control and optimal power flow settings; the development of concrete applications in this context is a subject of an ongoing work. We also note that the proposed approach may also be useful in the context of continuation analysis \cite{ContinuationBook,Continuation1,Continuation2}, which could be of future research interest. Lastly, the extension of our analysis approach to the case of active voltage regulators and capacitor banks is another future research direction. }  

%\vspace*{-2cm}
\appendix

\subsection{Proof of Lemma \ref{lem:norm}}

%\begin{proof}
We need to show the three norm axioms. Trivially, note that $\xi(a \bs) = |a| \xi(\bs)$ for any $a \in \comps$. Next, the triangle inequality holds because
\footnotesize
\begin{align*}
&\xi(\bs + \bs') = \left\|\bW^{-1} \bY_{LL}^{-1} \bW^{-1} \diag(\bs^Y + \bs'^Y) \right\|_{\infty} \\
&\quad + \left\|\bW^{-1} \bY_{LL}^{-1} \bH^\sfT \diag (\bL|\bw|)^{-1} \diag(\bs^\Delta + \bs'^\Delta) \right\|_{\infty} \\
& = \left\|\bW^{-1} \bY_{LL}^{-1} \bW^{-1} \diag(\bs^Y)  + \bW^{-1} \bY_{LL}^{-1} \bW^{-1}\diag(\bs'^Y) \right\|_{\infty} \\
&\quad + \Big\|\bW^{-1} \bY_{LL}^{-1} \bH^\sfT \diag (\bL|\bw|)^{-1} \diag(\bs^\Delta) + \\ 
&\quad \quad \quad \quad \quad  \bW^{-1} \bY_{LL}^{-1} \bH^\sfT \diag (\bL|\bw|)^{-1} \diag(\bs'^\Delta) \Big\|_{\infty}\\
&\leq \xi^Y(\bs) + \xi^Y(\bs') + \xi^\Delta(\bs)  + \xi^\Delta(\bs') = \xi(\bs) + \xi(\bs'),
\end{align*}
\normalsize
where the inequality follows by the triangle inequality for the induced matrix norm. Finally, if $\xi(\bs) = 0$, it necessarily holds that $\bW^{-1} \bY_{LL}^{-1} \bW^{-1} \diag(\bs^Y)$ and $\bW^{-1} \bY_{LL}^{-1} \bH^\sfT \diag (\bL|\bw|)^{-1} \diag(\bs^\Delta)$ are zero matrices. This necessarily implies that $\bs^Y$ and $\bs^\Delta$ are zero vectors.
%\end{proof}

\subsection{Proof of Theorem \ref{thm:exist_general}}

For the purpose of the proof, we find it convenient to re-parametrize using $\bu := \bW^{-1}\bv$. Then, \eqref{eqn:fp} is equivalent to

\begin{align} \label{eqn:lf_u}
\bu=  \tilde{\mathbf{G}}_{\bs^Y\bs^\Delta}(\bu) &=\bone + \bW^{-1} \bY_{LL}^{-1}\conj{\bW}^{-1} \diag(\conj{\bu})^{-1}\conj{\bs^{Y}} \nonumber\\
&\, + \bW^{-1} \bY_{LL}^{-1} \bH^\sfT \diag\left(\bH \conj{\bW} \conj{\bu} \right)^{-1} \conj{\bs^{\Delta}}.
\end{align}

As $\bW$ defines an invertible relationship between $\bv$ and $\bu$, we next focus on the solution properties of \eqref{eqn:lf_u}. By the Banach fixed-point theorem, what we need to show is that $\tilde{\mathbf{G}}_{\bs^Y\bs^\Delta}(\bu)$ is a self-mapping and contraction mapping on
\begin{equation} \label{eqn:D_tilde}
\tilde{\mathcal{D}}_\rho(\hat{\bu}) := \{\bu: \, |(\bu)_j - (\hat{\bu})_j| \leq \rho, \, j = 1 \ldots 3N \}
\end{equation}
for some $\rho\in(0,\gamma(\hat{\bv}))$ that satisfies \eqref{eqn:cond1_gen} and \eqref{eqn:cond2_gen}.

\subsubsection{Proof of Self-Mapping}
The goal here is to show that, for $\rho\in(0,\gamma(\hat{\bv}))$ fulfilling \eqref{eqn:cond1_gen}, $\|\bu^{(k)}-\hat{\bu}\|_\infty\leq\rho$ leads to $\|\bu^{(k+1)}-\hat{\bu}\|_\infty\leq\rho$.

By definition, we have
%\begin{footnotesize}
\footnotesize
\begin{align} \label{eqn:self_map}
&\mathbf{u}^{(k+1)}-\hat{\mathbf{u}} \nonumber \\
&=  \mathbf{W}^{-1}\mathbf{Y}_{LL}^{-1}\overline{\mathbf{W}}^{-1}\left(\diag(\overline{\mathbf{u}}^{(k)})^{-1}\overline{\mathbf{s}^Y}-\diag(\overline{\hat{\mathbf{u}}})^{-1}\overline{\hat{\mathbf{s}}^Y}\right) \nonumber \\
& +\mathbf{W}^{-1}\mathbf{Y}_{LL}^{-1}\mathbf{H}^T\left(\diag(\mathbf{H}\overline{\mathbf{W}}\overline{\mathbf{u}}^{(k)})^{-1}\overline{\mathbf{s}^\Delta}-\diag(\mathbf{H}\overline{\mathbf{W}}\overline{\hat{\mathbf{u}}})^{-1}\overline{\hat{\mathbf{s}}^\Delta}\right) \nonumber \\
&=  \mathbf{W}^{-1}\mathbf{Y}_{LL}^{-1}\overline{\mathbf{W}}^{-1}\left(\diag(\overline{\mathbf{u}}^{(k)})^{-1}\overline{\mathbf{s}^Y}-\diag(\overline{\mathbf{u}}^{(k)})^{-1}\overline{\hat{\mathbf{s}}^Y}\right) \nonumber \\
& +\mathbf{W}^{-1}\mathbf{Y}_{LL}^{-1}\overline{\mathbf{W}}^{-1}\left(\diag(\overline{\mathbf{u}}^{(k)})^{-1}\overline{\hat{\mathbf{s}}^Y}-\diag(\overline{\hat{\mathbf{u}}})^{-1}\overline{\hat{\mathbf{s}}^Y}\right) \nonumber \\
& +\mathbf{W}^{-1}\mathbf{Y}_{LL}^{-1}\mathbf{H}^T\left(\diag(\mathbf{H}\overline{\mathbf{W}}\overline{\mathbf{u}}^{(k)})^{-1}\overline{\mathbf{s}^\Delta}-\diag(\mathbf{H}\overline{\mathbf{W}}\overline{\mathbf{u}}^{(k)})^{-1}\overline{\hat{\mathbf{s}}^\Delta}\right) \nonumber \\
& +\mathbf{W}^{-1}\mathbf{Y}_{LL}^{-1}\mathbf{H}^T\left(\diag(\mathbf{H}\overline{\mathbf{W}}\overline{\mathbf{u}}^{(k)})^{-1}\overline{\hat{\mathbf{s}}^\Delta}-\diag(\mathbf{H}\overline{\mathbf{W}}\overline{\hat{\mathbf{u}}})^{-1}\overline{\hat{\mathbf{s}}^\Delta}\right).
\end{align}
%\end{footnotesize}
\normalsize
We can rearrange the right-hand side of \eqref{eqn:self_map} as follows.
For example, for the second term, we have
\iffalse
\begin{align} \label{eqn:triangA}
& \mathbf{W}^{-1}\mathbf{Y}_{LL}^{-1}\overline{\mathbf{W}}^{-1}\left(\diag(\overline{\mathbf{u}}^{(k)})^{-1}\overline{\mathbf{s}^Y}-\diag(\overline{\mathbf{u}}^{(k)})^{-1}\overline{\hat{\mathbf{s}}^Y}\right) \nonumber \\
= & \bW^{-1}\mathbf{Y}_{LL}^{-1}\overline{\bW}^{-1}\diag(\overline{\bs^Y}-\overline{\hat{\bs}^Y})
\begin{bmatrix}
\displaystyle\frac{1}{(\overline{\bu}^{(k)})_1} \\
\vdots \\
\displaystyle\frac{1}{(\overline{\bu}^{(k)})_{3N}}
\end{bmatrix}
,
\end{align}
and
\fi
\footnotesize
\begin{align} \label{eqn:triangB}
&\mathbf{W}^{-1}\mathbf{Y}_{LL}^{-1}\overline{\mathbf{W}}^{-1}\left(\diag(\overline{\mathbf{u}}^{(k)})^{-1}\overline{\hat{\mathbf{s}}^Y}-\diag(\overline{\hat{\mathbf{u}}})^{-1}\overline{\hat{\mathbf{s}}^Y}\right) \nonumber \\
&=-\bW^{-1}\mathbf{Y}_{LL}^{-1}\overline{\bW}^{-1}\diag(\overline{\hat{\bs}^Y})
\begin{bmatrix}
\displaystyle\frac{(\overline{\bu}^{(k)})_1-(\overline{\hat{\bu}})_1}{(\overline{\bu}^{(k)})_1(\overline{\hat{\bu}})_1}
\ldots 
\displaystyle\frac{(\overline{\bu}^{(k)})_{3N}-(\overline{\hat{\bu}})_{3N}}{(\overline{\bu}^{(k)})_{3N}(\overline{\hat{\bu}})_{3N}}
\end{bmatrix}^\sfT
\end{align}
\normalsize
Similar rearrangements can be applied to the remaining terms in \eqref{eqn:self_map}.
%It is straightforward to draw analogy for the other two lines in \eqref{eqn:self_map}.
Therefore, by triangular inequality, the definition of the induced matrix norm, and definition \eqref{eqn:xi}, it holds that
\footnotesize
\begin{align}
& \|\mathbf{u}^{(k+1)}-\hat{\mathbf{u}}\|_\infty 
\leq  \xi^Y(\bs-\hat{\bs})\max\limits_{j\in\{1,...,3N\}}\displaystyle\left (1/|(\bu^{(k)})_j| \right) \nonumber \\
& +\xi^Y(\hat{\bs})\max\limits_{j\in\{1,...,3N\}}\displaystyle \left( |(\bu^{(k)})_j-(\hat{\bu})_j|/[|(\bu^{(k)})_j||(\hat{\bu})_j|] \right)\nonumber \\
& +\xi^\Delta(\bs-\hat{\bs})\max\limits_{j\in\{1,...,3N\}}\displaystyle\left(1 \middle /\displaystyle\frac{|(\bH\bW\bu^{(k)})_j|}{(\bL|\bw|)_j} \right) \nonumber \\
& +\xi^\Delta(\hat{\bs})\max\limits_{j\in\{1,...,3N\}}\displaystyle\displaystyle\frac{|(\bH\bW\bu^{(k)})_j-(\bH\bW\hat{\bu})_j|/(\bL|\bw|)_j }{\displaystyle\frac{|(\bH\bW\bu^{(k)})_j|}{(\bL|\bw|)_j}\displaystyle\frac{|(\bH\bW\hat{\bu})_j|}{(\bL|\bw|)_j}}.
\end{align}
\normalsize
Observe that the following is true for any $j\in\{1,...,3N\}$ whenever $\|\bu^{(k)}-\hat{\bu}\|\leq\rho$:
\begin{subequations} \label{eqn:ineqAll}
\begin{align}
& |(\bu^{(k)})_j|\geq|(\hat{\bu})_j|-|(\bu^{(k)})_j-(\hat{\bu})_j|\geq\alpha(\hat{\bv})-\rho \label{eqn:ineqU} \\
& |(\bH\bW\bu^{(k)})_j-(\bH\bW\hat{\bu})_j| \leq (\bL|\bw|)_j\rho \label{eqn:ineqA} \\
& |(\bH\bW\bu^{(k)})_j|\geq(\beta(\hat{\bv})-\rho)(\bL|\bw|)_j, \label{eqn:ineqB}
\end{align}
\end{subequations}
where $\alpha(\cdot)$ and $\beta(\cdot)$ are defined in \eqref{eqn:gamma}.
In details, \eqref{eqn:ineqA} holds because
\footnotesize
\begin{align}
& |(\bH\bW\bu^{(k)})_j-(\bH\bW\hat{\bu})_j| \nonumber \\
&=  \left|\left((\bw)_\ell(\bu^{(k)})_\ell-(\bw)_{\ell'}(\bu^{(k)})_{\ell'}\right)-\left((\bw)_\ell(\hat{\bu})_\ell-(\bw)_{\ell'}(\hat{\bu})_{\ell'}\right)\right| \nonumber \\
&\leq |(\bw)_\ell|\left|(\bu^{(k)})_\ell-(\hat{\bu})_\ell\right|+|(\bw)_{\ell'}|\left|(\bu^{(k)})_{\ell'}-(\hat{\bu})_{\ell'}\right| \nonumber \\
&\leq  \left(|(\bw)_\ell|+|(\bw)_{\ell'}|\right)\|\bu^{(k)}-\hat{\bu}\|_\infty \leq  (\bL|\bw|)_j\rho
\end{align}
\normalsize
for some $\ell,\ell'$ in $\{1,...,3N\}$, and \eqref{eqn:ineqB} holds because
\begin{align}
& |(\bH\bW\bu^{(k)})_j|
\geq  |(\bH\bW\hat{\bu})_j|-|(\bH\bW\bu^{(k)})_j-(\bH\bW\hat{\bu})_j| \nonumber \\
&\geq  (\beta(\hat{\bv})-\rho)(\bL|\bw|)_j.
\end{align}
In this way, for $\rho\in(0,\gamma(\hat{\bv}))$, we obtain 
\begin{align}
& \|\mathbf{u}^{(k+1)}-\hat{\mathbf{u}}\|_\infty \nonumber \\
& \leq\frac{\xi^Y(\mathbf{s}-\hat{\mathbf{s}})+\displaystyle\rho\xi^Y(\hat{\mathbf{s}})/\alpha(\hat{\bv})}{\alpha(\hat{\bv})-\rho}+\frac{\xi^\Delta(\mathbf{s}-\hat{\mathbf{s}})+\displaystyle\rho\xi^\Delta(\hat{\mathbf{s}})/\beta(\hat{\bv})}{\beta(\hat{\bv})-\rho}.
\end{align}
This implies that $\|\bu^{(k)}-\hat{\bu}\|_\infty\leq\rho$ gives $\|\bu^{(k+1)}-\hat{\bu}\|_\infty\leq\rho$ for $\rho\in(0,\gamma(\hat{\bv}))$ fulfilling \eqref{eqn:cond1_gen}, and hence completes the proof.

\subsubsection{Proof of Contraction}
In this part, assuming there is a $\rho\in(0,\gamma(\hat{\bv}))$ fulfilling \eqref{eqn:cond1_gen}, we prove that $\|\bu^{(k+1)}-\bu^{(k)}\|_\infty<\|\bu^{(k)}-\bu^{(k-1)}\|_\infty$ if $\rho$ further satisfies \eqref{eqn:cond2_gen}.

Similar to the proof of self-mapping, we have
\begin{footnotesize}
\begin{align*}
& \mathbf{u}^{(k+1)}-\mathbf{u}^{(k)} \nonumber \\
&=  \mathbf{W}^{-1}\mathbf{Y}_{LL}^{-1}\overline{\mathbf{W}}^{-1}\left(\diag(\overline{\mathbf{u}}^{(k)})^{-1}\overline{\mathbf{s}^\mathbf{Y}}-\diag(\overline{\mathbf{u}}^{(k-1)})^{-1}\overline{\mathbf{s}^\mathbf{Y}}\right) \nonumber \\
& +\mathbf{W}^{-1}\mathbf{Y}_{LL}^{-1}\mathbf{H}^T\left(\diag(\mathbf{H}\overline{\mathbf{W}}\overline{\mathbf{u}}^{(k)})^{-1}\overline{\mathbf{s}^\Delta}-\diag(\mathbf{H}\overline{\mathbf{W}}\overline{\mathbf{u}}^{(k-1)})^{-1}\overline{\mathbf{s}^\Delta}\right).
\end{align*}
\end{footnotesize} \vspace*{-0.3cm}
Then, via derivations analogues to \eqref{eqn:triangB} and \eqref{eqn:ineqAll}, there is
\footnotesize
\begin{align} \label{eqn:contraction}
& \|\mathbf{u}^{(k+1)}-\mathbf{u}^{(k)}\|_\infty \nonumber \\
&\leq  \xi^Y(\bs)\max\limits_{j\in\{1,\ldots,3N\}}\displaystyle|(\bu^{(k)})_j-(\bu^{(k-1)})_j|/[|(\bu^{(k)})_j||(\bu^{(k-1)})_j|] \nonumber \\
& +\xi^\Delta(\bs)\max\limits_{j\in\{1,...,3N\}}\displaystyle\frac{\displaystyle|(\bH\bW\bu^{(k)})_j-(\bH\bW\bu^{(k-1)})_j|/(\bL|\bw|)_j}{\displaystyle\frac{|(\bH\bW\bu^{(k)})_j|}{(\bL|\bw|)_j}\displaystyle\frac{|(\bH\bW\bu^{(k-1)})_j|}{(\bL|\bw|)_j}} \nonumber \\
&\leq  \left(\frac{\xi^Y(\mathbf{s})}{(\alpha(\hat{\bv})-\rho)^2}+\frac{\xi^\Delta(\mathbf{s})}{(\beta(\hat{\bv})-\rho)^2}\right)\|\mathbf{u}^{(k)}-\mathbf{u}^{(k-1)}\|_\infty.
\end{align}
\normalsize
Clearly, $\|\bu^{(k+1)}-\bu^{(k)}\|_\infty<\|\bu^{(k)}-\bu^{(k-1)}\|_\infty$ if $\rho$ further satisfies \eqref{eqn:cond2_gen}.

\vspace*{-0.5cm}
\subsection{Proof of Theorem \ref{thm:exist_lf}}
%Like in the proof of Theorem \ref{thm:exist_general}, we rely on Banach Fixed Point Theorem. However, most of the reasoning here will be built on top of Theorem \ref{thm:exist_general}.

\revv{For item (i), we first note that because the Jacobian  $\bJ$ associated with the mapping $\bh$  in \eqref{eqn:real_lf} is a square matrix, the existence and uniqueness of the solution to the set of linear equations $\bJ \frac{\partial \by}{\partial \bx} = \bI_{12N}$ is equivalent to the invertibility of $\bJ$. In such case, the solution is given by $ \frac{\partial \by}{\partial \bx} = \bJ^{-1}$. Therefore, we can analyze the invertibility of $\bJ$ by analyzing the set of equations \eqref{eqn:part_der}.
In particular, we next show that if \eqref{eqn:cond1} is satisfied, \eqref{eqn:part_der} has a unique solution.
Because the system is linear with respect to the rectangular coordinates and there are as many unknowns as equations, the result is equivalent to showing that the corresponding homogeneous system of equations has only the trivial solution (see, e.g., \cite{H90}). 
Note that the homogeneous system is the same for every  column of \eqref{eqn:part_der} and is given by
\begin{subequations} \label{eqn:part_der_hom}
\begin{align} 
&\diag\left(\bH^\sfT \conj{\bi^{\Delta}}\right) \bDelta_V + \diag (\bv) \bH^\sfT \conj{\bDelta}_I,  \notag\\
& \quad = \diag(\bv) \conj{\bY}_{LL} \conj{\bDelta}_V + \diag (\conj{\bY}_{L0} \conj{\bv}_0 + \conj{\bY}_{LL} \conj{\bv}) \bDelta_V
 \\
&{\bf 0} = \diag\left(\bH \bv \right)\conj{\bDelta}_I + \diag(\conj{\bi^{\Delta}})\bH \bDelta_V,
\end{align}
\end{subequations}
where $\bDelta_V, \bDelta_I$ are solution vectors. 
	
Assume, by the way of contradiction, that there exists a solution $\bDelta' := ({\bDelta'}_V^{\sfT}, {\bDelta'}_I^{\sfT})^\sfT$ to \eqref{eqn:part_der_hom} such that $\bDelta' \neq 0$. In particular, any vector $\bDelta^{\epsilon} := \epsilon \bDelta'$ for $\epsilon > 0$ is a solution to \eqref{eqn:part_der_hom}.

Now consider two  power networks with the same topology but different voltages and between-phase currents.  In particular, let
%\begin{align*}
$\bv_1^\epsilon = \bv + \bDelta^\epsilon_V$, %\quad 
$\bi^{\epsilon, \Delta}_1 = \bi^\Delta + \bDelta_{I}^\epsilon$, %\\
$\bv_2^\epsilon = \bv - \bDelta^\epsilon_V$, % \quad 
and
$\bi^{\epsilon, \Delta}_2 = \bi^\Delta - \bDelta_{I}^\epsilon$,
%\end{align*}
while $\bv_0$ is the same in both networks. Note that there exists $\epsilon_1 > 0$ such that for all $\epsilon < \epsilon_1$, $\bv_1^\epsilon, \bv_2^\epsilon \in \cD_{\rho^\dagger}(\bv)$, where $\cD_{\rho^\dagger}$ is defined in \eqref{eqn:D} (with $\hat{\bv} = \bv$).

Let $\bs_1^{\epsilon, Y}, \bs_1^{\epsilon,\Delta}, \bs_2^{\epsilon,Y}, \bs_2^{\epsilon,\Delta}$ be the corresponding power injections.  Using \eqref{eqn:lf}, we obtain that
\begin{align*}
&\bs_1^{\epsilon,Y} - \bs_2^{\epsilon,Y} =  2\Big( \diag(\bv) \conj{\bY}_{LL} \conj{\bDelta}^\epsilon_V \\
& \quad + \diag (\conj{\bY}_{L0} \conj{\bv}_0 + \conj{\bY}_{LL} \conj{\bv}) \bDelta^\epsilon_V 
- \diag\left(\bH^\sfT \conj{\bi^{\Delta}}\right) \bDelta^\epsilon_V\\ 
& \quad - \diag (\bv) \bH^\sfT \conj{\bDelta}^\epsilon_I   \Big), \\
& \bs_1^{\epsilon, \Delta} - \bs_2^{\epsilon, \Delta} =  2\left(\diag\left(\bH \bv \right)\conj{\bDelta}^\epsilon_I + \diag(\conj{\bi^{\Delta}})\bH \bDelta^\epsilon_V \right),
\end{align*}
which by \eqref{eqn:part_der_hom} implies that
$
\bs_1^{\epsilon, Y} = \bs_2^{\epsilon, Y}$ and $  \bs_1^{\epsilon, \Delta} = \bs_2^{\epsilon, \Delta}.
$

Let $\bs^{\epsilon} := ((\bs_1^{\epsilon, Y})^\sfT, (\bs_1^{\epsilon, \Delta})^\sfT)^\sfT$. It is easy to see that there exists $\epsilon_2 > 0$ such that for all $\epsilon < \epsilon_2$, $\bs^\epsilon$ satisfies \eqref{eqn:cond2} (with $\hat{\bs} = \bs$). Let $\epsilon^* = \min \{\epsilon_1, \epsilon_2 \}$. Then, by Theorem \ref{thm:exist_lf}, we have that for any $\epsilon \in (0, \epsilon^*)$, $\bv^\epsilon_1 = \bv^\epsilon_2$ and $\bi^{\epsilon, \Delta}_1 = \bi^{\epsilon, \Delta}_2$. This is equivalent to having $\bDelta_V^\epsilon = {\bf 0}$ and $\bDelta_I^\epsilon = {\bf 0}$, which is a contradiction to our assumption that $\bDelta^\epsilon \neq 0$. This completes the proof of item (i).
}

\revv{For items (ii)-(iv) in the theorem,} we show that conditions \eqref{eqn:cond1} and \eqref{eqn:cond2} imply conditions \eqref{eqn:cond1_gen} and \eqref{eqn:cond2_gen} of Theorem \ref{thm:exist_general}. From the proof of Lemma 1 in \cite{congLF}, whenever \eqref{eqn:cond1} and \eqref{eqn:cond2} are satisfied, we have
\begin{equation}
\rho^2-\rho ((\gamma(\hat{\bv}))^2 - \xi(\hat{\bs}))/\gamma(\hat{\bv}) +\xi(\bs-\hat{\bs})\leq0
\end{equation}
for $\rho\in[\rho^\dagger(\hat{\bv},\hat{\bs},\bs),\rho^\ddagger(\hat{\bv},\hat{\bs})]\subseteq(0,\gamma(\hat{\bv}))$. %according to the results in \cite{congThreePhase}. 
After re-organization, the above inequality becomes
\begin{equation}
\displaystyle\frac{\xi(\bs-\hat{\bs})+\displaystyle\rho\xi(\hat{\bs})/\gamma(\hat{\bv})}{\gamma(\hat{\bv})-\rho}\leq\rho.
\end{equation}
Note that
\footnotesize
\begin{equation}
\frac{\xi^Y(\mathbf{s}-\hat{\mathbf{s}})+\displaystyle\frac{\xi^Y(\hat{\mathbf{s}})}{\alpha(\hat{\bv})}\rho}{\alpha(\hat{\bv})-\rho}+\frac{\xi^\Delta(\mathbf{s}-\hat{\mathbf{s}})+\displaystyle\frac{\xi^\Delta(\hat{\mathbf{s}})}{\beta(\hat{\bv})}\rho}{\beta(\hat{\bv})-\rho}\leq\displaystyle\frac{\xi(\bs-\hat{\bs})+\displaystyle\frac{\xi(\hat{\bs})}{\gamma(\hat{\bv})}\rho}{\gamma(\hat{\bv})-\rho},
\end{equation}
\normalsize
and hence \eqref{eqn:cond1_gen} is satisfied. Namely, $\tilde{\mathbf{G}}_{\bs^Y\bs^\Delta}(\bu)$ is a self-mapping on domain $\tilde{\mathcal{D}}_\rho(\hat{\bu})$ for $\rho\in[\rho^\dagger(\hat{\bv},\hat{\bs},\bs),\rho^\ddagger(\hat{\bv},\hat{\bs})]$.

Further, taking into account \eqref{eqn:cond2} and the fact that $\xi(\cdot)$ is a norm (cf.~Lemma \ref{lem:norm}), we have
\footnotesize
\begin{align} \label{eqn:Nonsing}
\xi(\bs) \nonumber &\leq  \xi(\bs-\hat{\bs})+\xi(\hat{\bs}) 
<  \frac{1}{4}\left(\frac{(\gamma(\hat{\bv}))^2 - \xi(\hat{\bs})}{\gamma(\hat{\bv})} \right)^2+\xi(\hat{\bs}) \nonumber \\
&=  \frac{1}{4}\left(\frac{(\gamma(\hat{\bv}))^2 + \xi(\hat{\bs})}{\gamma(\hat{\bv})} \right)^2 =  (\gamma(\hat{\bv})-\rho^\ddagger(\hat{\bv},\hat{\bs}))^2.
\end{align}
\normalsize
Then, it is easy to see that
\footnotesize
\[
\frac{\xi^Y(\mathbf{s})}{(\alpha(\hat{\bv})-\rho)^2}+\frac{\xi^\Delta(\mathbf{s})}{(\beta(\hat{\bv})-\rho)^2} 
\leq  \displaystyle\frac{\xi(\bs)}{\left(\gamma(\hat{\bv})-\rho\right)^2} <  \displaystyle\frac{\left(\gamma(\hat{\bv})-\rho^\ddagger(\hat{\bv},\hat{\bs})\right)^2}{\left(\gamma(\hat{\bv})-\rho\right)^2} \leq 1
\]
\normalsize
for $\rho\in[\rho^\dagger(\hat{\bv},\hat{\bs},\bs),\rho^\ddagger(\hat{\bv},\hat{\bs})]$. It implies that $\tilde{\mathbf{G}}_{\bs^Y\bs^\Delta}(\bu)$ is also a contraction mapping on $\tilde{\mathcal{D}}_\rho(\hat{\bu})$ with $\rho\in[\rho^\dagger(\hat{\bv},\hat{\bs},\bs),\rho^\ddagger(\hat{\bv},\hat{\bs})]$.
This completes the proof \revv{of items (ii)-(iv)}.

\revv{For item (v)}, we derive as follows. First, by \eqref{eqn:gamma} and the bound of \eqref{eqn:ineqU}, there is
\begin{equation} \label{eqn:NonsingA}
\gamma(\hat{\bv})-\rho^\ddagger(\hat{\bv},\hat{\bs})\leq\alpha(\hat{\bv})-\rho^\ddagger(\hat{\bv},\hat{\bs})\leq\alpha(\bv).
\end{equation}
Then, by the bound of \eqref{eqn:ineqB}, we have
\begin{equation} \label{eqn:NonsingB}
\gamma(\hat{\bv})-\rho^\ddagger(\hat{\bv},\hat{\bs})\leq\beta(\hat{\bv})-\rho^\ddagger(\hat{\bv},\hat{\bs})\leq\beta(\bv).
\end{equation}
Combination of \eqref{eqn:NonsingA} and \eqref{eqn:NonsingB} yields
\begin{equation} 
\gamma(\hat{\bv})-\rho^\ddagger(\hat{\bv},\hat{\bs})\leq\gamma(\bv).
\end{equation}
Therefore, by \eqref{eqn:Nonsing}, we have
\begin{equation}
\xi(\bs)<(\gamma(\hat{\bv})-\rho^\ddagger(\hat{\bv},\hat{\bs}))^2\leq(\gamma(\bv))^2
\end{equation}
and hereby complete the proof.

\subsection{Proof of Theorem \ref{thm:fixedAccuracy}}
%\begin{proof}
Note that \eqref{eqn:lin_fp} is in fact a single iteration of the fixed-point equation initialized at $\hat{\bv}$. 
Therefore, by identifying $\bv^{(0)} = \hat{\bv}$ and $\bv^{(1)} = \tilde{\bv}$, we have that
\begin{align*}
\|\tilde{\bv} - \bv \|_\infty \leq q \|\hat{\bv} - \bv \|_\infty &\leq q \| \bw \|_\infty \|\hat{\bu} - \bu \|_\infty \leq q \| \bw \|_\infty \rho^\dagger,
\end{align*}
where $q < 1$ is the contraction coefficient given in the proof of Theorem \ref{thm:exist_lf} -- cf.~\eqref{eqn:contraction}; the first inequality follows by the Banach fixed point theorem; the second inequality holds  by definition of $\bv = \bW \bu$; and the last inequality follows because $\bv \in \cD_{\rho^\dagger}(\hat{\bv})$ (cf.~\eqref{eqn:D}). 
%\end{proof}

%%%%%%%%%%%%%%%%%%%%%%%%%%%%%%%%%%%%%%%%%%%%%
\bibliographystyle{IEEEtran}
\bibliography{biblio.bib}
\end{document}